\documentclass[reqno, 12pt]{amsart}

\usepackage{amssymb,amsmath}
\usepackage{amsthm}
\usepackage{upref}
\usepackage{enumerate}
\usepackage{amsfonts,bbold}

 \makeatletter
\def\LaTeX{\leavevmode L\raise.42ex
    \hbox{\kern-.3em\size{\sf@size}{0pt}\selectfont A}\kern-.15em\TeX}
 \makeatother
 
\DeclareMathOperator{\clos}{clos}

\numberwithin{equation}{section}

\newtheorem{lemma}{Lemma}[section]
\newtheorem{theorem}[lemma]{Theorem} 
\newtheorem{corollary}[lemma]{Corollary}
\newtheorem{proposition}[lemma]{Proposition}

\theoremstyle{definition}

\newtheorem{definition}[lemma]{Definition}
\newtheorem{example}[lemma]{Example}

\newtheorem{remark}[lemma]{Remark}

\newcommand{\dist}{\operatorname{dist}}

\renewcommand{\det}{\operatorname{Det}}
\newcommand{\tr}{\operatorname{Tr}}

 \DeclareMathOperator*{\slim}{s-lim}

 \newcommand{\supp}{\operatorname{supp}}
  \newcommand{\e}{\eqref}
\newcommand{\q}{\quad}

\renewcommand{\d}{\delta}

\newcommand{\ov}{\overline}
\newcommand{\z}{\zeta}
\newcommand{\wt}{\widetilde}

   \newcommand{\var}{\operatorname{var}}

\renewcommand\Im{\operatorname{Im}}
\renewcommand\Re{\operatorname{Re}}

\newenvironment{pf}{\begin{proof}}{\end{proof}}

\def\qqq{\mathrel{\subset\mkern-15mu\lower.38ex\hbox{${\scriptscriptstyle\rightarrow}$}}}

\let\cal\mathcal

\let\Bbb\mathbb

\begin{document}
\title[Scattering theory   for Jacobi operators]
 {Analytic scattering theory   for Jacobi operators and Bernstein-Szeg\"o asymptotics of orthogonal polynomials}
\author{ D. R. Yafaev  }


\dedicatory{To the memory of Lyudvig Dmitrievich Faddeev}

\address{ IRMAR, Universit\'{e} de Rennes I, Campus de
  Beaulieu,  Rennes, 35042  FRANCE and SPGU, Univ. Nab. 7/9, Saint Petersburg, 199034 RUSSIA}
\email{yafaev@univ-rennes1.fr}
\subjclass[2000]{33C45, 39A70,  47A40, 47B39}

  \keywords {Jacobi matrices, the discrete Schr\"odinger operator,  orthogonal polynomials, asymptotics for large numbers, the Szeg\"o function}

\thanks {Supported by  project   Russian Science Foundation   17-11-01126}

\begin{abstract}
We study  semi-infinite Jacobi matrices $H=H_{0}+V$ corresponding to trace class perturbations $V$ of the ``free"
 discrete Schr\"odinger operator $H_{0}$. Our goal is to construct   various
 spectral quantities of the operator $H$, such as the weight function, eigenfunctions  of its continuous spectrum, the wave operators  for the pair $H_{0}$, $H$, the   scattering matrix, the spectral shift function, etc.
 This allows us to find the asymptotic behavior   of the orthonormal polynomials $P_{n}(z)$ associated to the Jacobi matrix $H $ as $n\to\infty$. In particular, we consider the case of $z$ inside the spectrum $[-1,1]$  of $H_{0}$  when this asymptotics has an oscillating character of the Bernstein-Szeg\"o type and the case of $z$ at the end points $\pm 1$.
   \end{abstract}

\maketitle

\thispagestyle{empty}

\section{Introduction}

{\bf 1.1.}
The theory of the Schr\"odinger operator $D^2   +b(x)$, $D=-i d/dx$,  with a short-range potential $b(x)$ is to a large extent due to L.~D.~Faddeev. We note, in particular, his classical papers \cite{F1, Finv} on the direct and inverse quantum scattering problems (both in the one- and multi-dimensional cases) and 
\cite{Fadd1} on the trace formulas for the   operator $D^2   +b(x)$. Later, the paper \cite{Fadd1} was significantly generalized by him jointly with V.~S.~Buslaev in \cite{BF}. Some of these results were exposed in the book  \cite{F-T} by L.~D.~Faddeev and L.~A.~Takhtajan.

 As is well known, the theories of Jacobi operators given by three-diagonal matrices   and of differential operators $D a (x)D +b(x)$   are  to a large extent similar. This is true   for Jacobi operators acting in the space  $\ell^2 ({\Bbb Z})$ and   differential operators
acting in the space  $L^2 ({\Bbb R})$ as well as for the corresponding operators acting in the spaces  $\ell^2 ({\Bbb Z}_{+})$ and    $L^2 ({\Bbb R}_{+})$, respectively. We  refer to the book \cite{Teschl}  where this analogy is described in a sufficiently detailed way.
Both classes of the operators are very important in physical applications.  Moreover, Jacobi operators   in the space  $\ell^2 ({\Bbb Z}_{+})$ are intimately related  (see, e.g.,  the classical book \cite{AKH}) to the theory of orthogonal polynomials. Necessary information on  orthogonal polynomials can be found in the other classical book \cite{Sz};  the  theory of Jacobi operators is carefully presented in  the  books \cite{AKH, Ber}.

Our goal is to develop the spectral and scattering theory for  Jacobi operators   in the space  $\ell^2 ({\Bbb Z}_{+})$ using basically the same approach as for the
differential operator $D a (x)D +b(x)$  in the space  $L^2 ({\Bbb R}_{+})$. This leads to new results for orthogonal polynomials.

\medskip

{\bf 1.2.}
  Jacobi operators    are defined in the space $\ell^2 ({\Bbb Z}_{+})$ by matrices
\begin{equation}
H= 
\begin{pmatrix}
 b_{0}&a_{0}& 0&0&0&\cdots \\
 a_{0}&b_{1}&a_{1}&0&0&\cdots \\
  0&a_{1}&b_{2}&a_{2}&0&\cdots \\
  0&0&a_{2}&b_{3}&a_{3}&\cdots \\
  \vdots&\vdots&\vdots&\ddots&\ddots&\ddots
\end{pmatrix} .
\label{eq:ZP+}\end{equation}
The entries $a_{n}$ here are arbitrary positive numbers, and $b_{n}$  are arbitrary real numbers. We denote by $e_{n}$, $ n\in {\Bbb Z}_{+}$,  the canonical  basis in $\ell^2 ({\Bbb Z}_{+})$, that is, all components of the vector $e_{n}$ are zeros, except the $n$-th component which equals $1$. If the sequences 
$\{a_{n}\}$ and $\{b_{n}\}$  are bounded, then the Jacobi operator $H$ is bounded   in $\ell^2 ({\Bbb Z}_{+})$. Let us denote by $dE(\lambda)$  the  spectral family of the self-adjoint operator $H$  and define the corresponding  spectral measure $d\rho(\lambda)=d(E(\lambda)e_{0}, e_{0})$.

Given matrix \e{eq:ZP+}, one constructs
 polynomials $P_{n } (z)$   by the recurrent relation
 \begin{equation}
 a_{n-1} P_{n-1} (z) +b_{n} P_{n } (z) + a_{n} P_{n+1} (z)= z P_n (z), \q n\in{\Bbb Z}_{+}, 
\label{eq:Py}\end{equation}
and the boundary conditions $P_{-1 } (z) =0$, $P_0 (z) =1$. Then $P_{n } (z)$ is a polynomial of degree $n$, that is,
 \begin{equation}
P_{n}(z)= k_{n}( z^{n}+ r_{n} z^{n-1}+\cdots).
\label{eq:norm}\end{equation}
  Comparing the coefficients at $z^{n+1}$ and $z^n$ in the left- and right-hand sides of \e{eq:Py}, we see that
  $k_{n}=(a_{0}a_{1}\cdots a_{n-1})^{-1}>0$ and
   \begin{equation}
a_{n} =\frac{ k_{n}}{ k_{n+1}}, \q b_{n} =r_{n}- r_{n+1}   .
\label{eq:norm1}\end{equation}
Obviously, $P(z)=\{ P_{n} (z)\}_{n=0}^\infty$ satisfies the equation $HP(z)=zP(z)$, that is, it is an ``eigenvector"  of the operator $H$.

  Putting together \e{eq:ZP+} and \e{eq:Py}, we find that
     \begin{equation}
e_{n}=P_{n} (H) e_{0}
\label{eq:J2}\end{equation}
 for all $n\in{\Bbb Z}_{+}$. It follows that
    \[
    d(E(\lambda)e_{n}, e_{m})= P_{n}(\lambda) P_{m} (\lambda) d\rho(\lambda),
    \]
    whence
     \begin{equation}
\int_{-\infty}^\infty P_{n}(\lambda) P_{m}(\lambda) d\rho(\lambda) =\d_{n,m};
\label{eq:J5}\end{equation}
as usual, $\d_{n,n}=1$ and $\d_{n,m}=0$ for $n\neq m$. 
Thus,  the  polynomials  $P_{n}(\lambda)$ are orthogonal and normalized  in the space $L^2 ({\Bbb R};d\rho)$. Formula \e{eq:J2} also shows that the set of all vectors $H^n e_{0}$, $n\in{\Bbb Z}_{+}$, is dense  in the space $\ell^2 ({\Bbb Z}_{+})$, and hence the spectrum  of the Jacobi operator $H$   is simple   with the generating vector $e_{0}$.

Alternatively,   $\{P_0 (\lambda),P_1 (\lambda),\ldots,P_{n} (\lambda),\ldots\}$ can be obtained by the Gram-Schmidt orthogonalization of the monomials $\{1,\lambda,\ldots,\lambda^n,\ldots\}$ in the space $L^2({\Bbb R}_{+}; d\rho)$; one also has to additionally require that $k_{n} >0$ in \e{eq:norm}.  This fact is known as the Favard theorem. 
In contract to the continuous case (see the paper \cite{Ge-Le} where an integral equation was used), the inverse problem of reconstructing the Jacobi operator $H$ given its spectral measure $d\rho(\lambda)$ admits a quite explicit solution. Indeed, let a measure $d\rho(\lambda)$
have a bounded and infinite support, and let $P_{n}(z)$ be the corresponding polynomials satisfying    \e{eq:norm} and \e{eq:J5}. Then the coefficients of the operator $H$ can be recovered by formulas \e{eq:norm1}.

\medskip

{\bf 1.3.}
In the particular case $a_{n}=1/2$ and $b_{n}=0$, the operator \e{eq:ZP+} is denoted $H_{0}$.
It plays the  role of the ``free"  differential operator $D^2$ in the space $L^2 ({\Bbb R}_{+})$ with the boundary condition $u(0)=0$. The operator $H_{0}$ can   be diagonalized explicitly. 

In this paper, we suppose that the perturbation $V=  H - H_{0}$  satisfies a ``short-range" assumption
 \begin{equation}
\sum_{n=0}^\infty (|  a_{n}-1/2|+|b_{n}|)<\infty. 
\label{eq:Tr}\end{equation} 
Then $V  $ belongs to the trace class $\mathfrak{S}_{1}$, and its trace norm   is equivalent to the sum \e{eq:Tr}.
Under assumption \e{eq:Tr} the   spectrum $\sigma(H)$ of the operator $H$ is absolutely continuous on the interval $(-1,1)$, but 
the operator $H$ may have   discrete spectrum $\sigma_{d}(H)$ (possibly, infinite) in ${\Bbb R}\setminus [-1,1]$; moreover,
the points $+1$  or $-1$ may be its eigenvalues.

We construct   various
 spectral quantities of the operator $H$, such as  the weight function  $w(\lambda)$ defined by the equation
    $d\rho(\lambda)=w(\lambda) d\lambda$ for $\lambda\in(-1,1)$, eigenfunctions  of its continuous spectrum, the perturbation determinant $\Delta(z)$  and
  the wave operators $W_{\pm} (H,H_{0}) $ for the pair $H_{0}$, $H$, the   scattering matrix
 $S(\lambda)$, the spectral shift function  $\xi(\lambda)$, etc.
Under assumption  \e{eq:Tr} the functions $w(\lambda)$,
  $S(\lambda)$ and $\xi(\lambda)$ are continuous in $\lambda\in(-1,1)$. Moreover, $w(\lambda)\neq 0$ and, as shown in \cite{KS}, the so-called Szeg\"o condition
    \begin{equation}
\int_{-1}^1 \ln w(\lambda) (1-\lambda^2)^{-1/2} d\lambda> -\infty 
\label{eq:Szeg}\end{equation}
is satisfied.

 Spectral results on Jacobi operators lead to the corresponding assertions for the polynomials $P_{n} (z)$ defined by \e{eq:Py}. For example, we show (see Theorem~\ref{Sz})  that, for $\lambda \in (-1,1)$,
 \begin{equation}
 P_{n} (\lambda)=  \sqrt{\frac{2}{\pi}} w(\lambda)^{-1/2}(1-\lambda^2)^{-1/4}  \sin ( (n+1)\arccos\lambda +\pi \xi(\lambda) ) + o(1) 
\label{eq:zF}\end{equation}
as $n\to\infty$. This asymptotic relation is the classical result of S.~Bernstein \cite{Bern} (see also     formula (12.1.8) in the G.~Szeg\"o book~\cite{Sz}). 
It is required in   \cite{Bern, Sz} that $\supp\rho\subset [-1,1]$ and   Lipschitz-Dini conditions are imposed on the weight function $w(\lambda)$. Under assumption \e{eq:Tr}  formula \e{eq:zF} is probably new. In  particular, we do no assume that  $\supp\rho\subset [-1,1]$. Note that the spectral shift function $\xi(\lambda)$ in formula \e{eq:zF}  is usually replaced in the theory of orthogonal polynomials by the so-called Szego function. 

As shows the example of Pollaczek polynomials (see formula \e{eq:Poll1}, below), assumption \e{eq:Tr}  cannot be significantly relaxed. For Pollaczek polynomials,  $ a_{n}-1/2  $ and $b_{n}  $ have order $n^{-1}$ as $n\to\infty$,
and the phase in formula \e{eq:zF} is essentially changed (see Section~5 in the Appendix to  \cite{Sz}). This resembles  the modification of the phase function for the Schr\"odinger operator with the Coulomb potential (see, e.g., formula (36,23)  in the book~\cite{LL}).

 As shown in the paper \cite{BF}  by V.~S.~Buslaev and L.~D.~Faddeev   (see Section~4.6 in \cite{YA}, for a detailed presentation), the expansion as $|z|\to\infty$ of the perturbation determinant $\Delta(z)$ for a pair $H_{0}$, $H$ of Schr\"odinger   operators contains  both integer and half-integer powers of $z^{-1}$. This leads to two series of trace identities: of integer and of  half-integer orders.  The first of them is stated in terms of integer powers of eigenvalues      and moments of the spectral shift function. The second series is stated in terms of half-integer powers of eigenvalues    and moments of the modulus of the perturbation determinant. 
For Jacobi operators  $H_{0}$, $H$,  the expansion of the perturbation determinant $\Delta(z)$  as $|z|\to\infty$ contains integer powers of $z^{-1}$ only. So there are no identities of half-integer orders. On the other hand, there is a version of    trace identities of integer order (known as Case sum rules)  stated in terms of the weight function, that is,  of the modulus of the perturbation determinant. We note that the Case sum rules involve
     Chebyshev polynomials of the operators $H_{0}$ and $H$ while the identities in spirit of \cite{BF} yield expressions for the traces $\tr(H^{n}-H_{0}^{n})$.
 
  This paper contains relatively few new results. However, we hope that a consistent  analogy between Jacobi and differential operators might shed a new light on  some aspects of the theory of orthogonal polynomials. A similar point of view was adopted in the paper \cite{KS}  devoted to Hilbert-Schmidt perturbations $V$ of the operator $H_{0}$.
  
   Note that,
in   \cite{KS}, necessary and sufficient conditions in terms of spectral data of $H$ were found for $V$ to be in the   Hilbert-Schmidt class. The problem of characterization of spectral data was also solved in
\cite{Ryck} for perturbations $V$  in the class (it is known now as the Ryckman class) relatively close to the trace  class $ \mathfrak{S}_{1}$.
It looks   tempting to obtain exhaustive results of such type for $V\in\mathfrak{S}_{1}$.

 \medskip

{\bf 1.4.}
The paper is organized as follows. In Section~2,  we consider equation \e{eq:Py} and, in addition to its polynomial solutions $P_{n}(z)$, we
introduce   so-called Jost solutions $f_{n}(z)$ of this equation. The Jost function $\Omega (z)$ basically coincides with $f_{-1}(z)$. The solutions $f_{n}(z)$  exponentially decay as $n\to\infty$ for $z\not\in\sigma (H_{0})=[-1,1]$ and oscillate for $z =\lambda\pm i0$ where $\lambda\in (-1,1)$.
Then a link between $f_{n}(\lambda\pm i0)$ and 
$P_{n}(\lambda)$ leads to asymptotic formulas for polynomials $P_{n}(\lambda)$. Thus,   we obtain formula \e{eq:zF}  (with the argument of $\Omega(\lambda\pm i0)$
playing the role of $\xi(\lambda)$); see Theorem~\ref{Sz}. Polynomial solutions $P_{n}(z)$ and  their asymptotics as $n\to\infty$ for $z\not\in [-1,1]$ are studied in Section~3.
Our proofs of the existence of $f_{n}(z)$ and uniform bounds on $P_{n}(z)$ rely on discrete ``Volterra integral" equations that are studied in the Appendix. 

Section~4 is devoted to  an investigation of the Jost solutions $f_{n}(z)$ and the polynomials $P_{n}(z)$ as $z\to \pm1$. This corresponds to the low energy scattering for the Schr\"odinger equation. In particular, we exhibit here examples of Jacobi operators with eigenvalues at the points $1$ or $-1$.

In Section~5, we construct the perturbation determinant $\Delta(z)$ and the spectral shift function for the pair $H_{0}$, $H$
and derive formulas for the traces $\tr(H^{n}-H_{0}^{n})$. We also show in Theorem~\ref{RRY} that the Jost function $\Omega (z)$ differs  from  $\Delta(z)$ by a numerical factor only. Scattering theory for  the pair $H_{0}$, $H$  is developed in Section~6.

A link of the Szeg\"o function with the perturbation determinant (or the Jost function)  is established in Section~7.
We here also show (again on the example of Pollaczek polynomials) that the weight function $w(\lambda)$ may tend to zero exponentially if assumption \e{eq:Tr} is slightly relaxed. Therefore the    Szeg\"o condition \e{eq:Szeg}   is sharp for trace class perturbations. In Section~8, the  results of the preceding sections are illustrated on the example of the Jacobi polynomials   when $ a_{n}-1/2  $ and $b_{n}  $ have order $n^{-2}$ as $n\to\infty$

Some parts  of this paper have non-trivial intersections with Section~2 of the paper \cite{KS} where, however, specific features of the Jacobi operators were extensively used. The analogy with the continuous case exploited in the present paper allows us to obtain some results for free.
 
 The modern approach based on a Riemann-Hilbert problem for matrix valued functions of \cite{F-I-K} combined with the steepest descent method of \cite{D-Z} is out of the scope of the present paper.


 \section{Jost solutions of the Jacobi equations and Bernstein-Szeg\"o asymptotics}

 {\bf 2.1.}
 In the canonical  basis
 $e_{n}$, $ n\in {\Bbb Z}_{+}$,   the Jacobi operator $H$ is 
  defined by the formula
 \begin{equation}
  H e_n= a_{n-1} e_{n-1} +b_{n} e_{n } + a_{n} e_{n+1}  , \q n\in {\Bbb Z}_{+} ,
\label{eq:J}\end{equation}
where we accept that $e_{-1}=0$. Relations \e{eq:ZP+} and \e{eq:J} are of course equivalent. The sequences $a_{n}>0$ and $b_{n}=\bar{b}_{n}$ where $n=0,1,\ldots$ are    assumed to be   bounded, so that $H$ is a bounded self-adjoint operator in   the space $\ell^2 ({\Bbb Z}_{+})$.

   Consider  now  the equation $H u =z u$, that is
 \begin{equation}
 a_{n-1} u_{n-1}  +b_{n} u_{n} + a_{n} u_{n+1}= z u_{n}, \q n\in{\Bbb Z}_{+}, 
\label{eq:Jy}\end{equation}
 for a sequence $u=\{u_{n} \}_{n=-1}^\infty$. The number $a_{-1}\neq 0$ can be chosen at our convenience; for definiteness, we put $a_{-1}=1/2$. Obviously, the values of $u_{k-1}$ and $u_{k }$ for some $k\in{\Bbb Z}_{+}$ determine the whole sequence $u_{n}$ satisfying equation \e{eq:Jy}.

 Let $f=\{ f_{n} \}_{n=-1}^\infty$ and $g=\{g_{n} \}_{n=-1}^\infty$ be two solutions of equation \e{eq:Jy}. A direct calculation shows that their Wronskian
  \begin{equation}
\{ f,g \}: = a_{n}  (f_{n}  g_{n+1}-f_{n+1}  g_{n})
\label{eq:Wr}\end{equation}
does not depend on $n=-1,0, 1,\ldots$. In particular, for $n=-1$ and $n=0$, we have
 \begin{equation}
\{ f,g \} = 2^{-1} (f_{-1}  g_{0}-f_{0}  g_{-1}) \q {\rm and} \q \{ f,g \} = a_{0}  (f_{0}  g_{ 1}-f_{ 1}  g_{0}).
\label{eq:Wr1}\end{equation}
Calculating the Wronskian \e{eq:Wr} for $n\to\infty$, we see that equation \e{eq:Jy} may have at most one (up to a multiplicative constant)  solution $u_{n}$ such that $u_{n}\to 0$ as $n\to\infty$.

  In the case $a_{n}=1/2$, $b_{n}=0$, the operator \e{eq:J} is known as the free discrete Schr\"odinger operator; it will be denoted $H_{0}$. The spectrum of the operator $H_{0}$ is simple, absolutely continuous and coincides with the interval $[-1,1]$. The corresponding spectral measure $d\rho_{0} (\lambda)= d(E_{0} (\lambda) e_{0}, e_{0})$ is given by the formula
 \[
d\rho_{0} (\lambda)= 2 \pi^{-1}  \sqrt{1-\lambda^2 }d\lambda,\q \lambda\in (-1,1). 
\]

Below we  fix the branch of the analytic function $\sqrt{z^2 -1}$ of $z\in {\Bbb C}\setminus [-1,1]=: \Pi$ by the condition $\sqrt{z^2 -1}>0$ for $z>1$. Obviously, this function is continuous on the closure $\clos{\Pi}$ of $\Pi$,    it equals $\pm i\sqrt{1-\lambda^2}$ for $z=\lambda\pm i0$, $\lambda\in (-1,1)$, and $\sqrt{z^2 -1}< 0$ for $z< -1$. Put
\begin{equation}
 \z (z) =z-\sqrt{z^2 -1}= (z + \sqrt{z^2 -1})^{-1};
 \label{eq:ome}\end{equation}
then $| \z (z)| <1$ for $z\in \Pi$.  
Since
\begin{equation}
2z= \z (z) +\z (z)^{-1} ,
\label{eq:omex}\end{equation}
the sequence $\{ \z (z)^n\}_{n=-1}^\infty$  satisfies the ``free" equation \e{eq:Jy}:
$ \z (z)^{n-1}+   \z (z)^{n+1}=2 z \z (z)^n$.
  For $\lambda\in [-1,1]$, it is common to set $\lambda=\cos \theta$ where $\theta\in [0,\pi ]$. Then 
\begin{equation}
\z (\lambda\pm i0)=e^{\mp i \theta}.
\label{eq:zt}\end{equation}

It is well known (see, e.g., \cite{Yjmp} for a detailed proof) that the  matrix elements of the resolvent $R_0 (z) =(H_0-z I)^{-1}$ (here and below $I$ is the identity operator)  of the operator $H_{0}$ are given by the formula
 \[
(R_{0}(z)e_{n}, e_{m})=  \frac{\z(z)^{n+m+2}-\z(z)^{|n-m|} }{\sqrt{z^2 -1}} 
\]
 for $z\in \Pi$ and all $n,m\in{\Bbb Z}_{+}$. In  particular, 
$(R_{0}(z)e_0, e_0)= - 2 \z(z) $.
The polynomials $P_{n}^{(0)} (z)$ defined by the recurrent equation \e{eq:Py} where $a_{n}=1/2$, $b_{n}=0$ and obeying the conditions $P_{-1}^{(0)} (z)=0$, $P_0^{(0)} (z)=1$  are   normalized  Chebyshev polynomials   of the second kind.  They satisfy the equation
 \begin{equation}
P_{n}^{(0)} (z)= \frac{1}{2\sqrt{z^2-1}} \big(\z(z)^{-n-1} -\z(z)^{n+1}\big)
\label{eq:OPR}\end{equation}
and, in particular,
 \begin{equation}
P_{n}^{(0)} (\lambda)= \frac{\sin((n+1)\theta)}{\sqrt{1-\lambda^2}}  ,\q \lambda=\cos\theta\in (-1,1).
\label{eq:OPs}\end{equation}

 Evidently, the ``perturbation" $V =H-H_{0}$ is given by the equality
\begin{equation}
  V e_n= (a_{n-1}-1/2) e_{n-1} +b_{n} e_{n } + (a_{n}-1/2) e_{n+1}   .
\label{eq:JV}\end{equation}
 If $a_{n}\to 1/2$ and $b_{n}\to 0$ as $n\to\infty$, then the operator $V$ is compact.
In this case the essential spectrum of the operator $H$ coincides with the interval $[-1,1]$, and its discrete spectrum $\sigma_{d} (H)$ consists of simple eigenvalues accumulating, possibly, to the points $1$ and $-1$ only.

Under assumption  \e{eq:Tr} the equation  \e{eq:Jy} has
  the so called Jost solution $ f(z)=\{ f_{n}(z)\}_{n=-1}^\infty$   distinguished by its  asymptotics as $n\to\infty$. 
  Our proof of this fact is   similar to the continuous case, but we give it for the completeness of our presentation in Appendix~A.1. It is based on the discrete Volterra integral equation
    \begin{equation}
  f_{n}(z)=\z(z)^n -\frac{1}{\sqrt{ z^2-1}}\sum_{m=n+1}^\infty (\z(z)^{n-m}- \z(z)^{m-n}) (Vf(z))_{m}
\label{eq:V}\end{equation}
where $z\in\Pi$  and $V$ is defined by formula \e{eq:JV}. 
A somewhat different proof can be found in \cite{KS} (see also the book \cite{Teschl}).  
   Put
  \begin{equation}
  \rho_{n}=\sum_{m=n}^\infty ( |a_m-1/2| + |b_m| ).
\label{eq:sig}\end{equation}
Then $ \rho_{n}\to 0$ as $n\to\infty$.

\begin{theorem}\label{Jost}
Let  assumption   \e{eq:Tr} be satisfied, and let $z\in\clos{\Pi}$, $z\neq \pm 1$. Then   equation \e{eq:Jy} has a   solution satisfying the condition
  \begin{equation}
  f_{n}(z)=\z(z)^n (1+O(\rho_{n}))
\label{eq:Jost}\end{equation}
as $n\to\infty$. Every function $f_{n}(z)$, $n=-1,0,1,\ldots$, depends analytically on $z\in\Pi$, and it is continuous in $z$ up to the cut along $[-1,1]$ except, possibly, the points $\pm 1$.
\end{theorem}

  Since $|\z (z)| <1$, it follows from \e{eq:Jost} that $ f_{n}(z)\to 0$ exponentially as $n\to \infty$ for $z\in \Pi$. In particular,  equation \e{eq:Jy} has only one   solution satisfying \e{eq:Jost}.
  Note that for the operator $H_{0}$, the error term in \e{eq:Jost} disappears and   the Jost function is $f_{n}(z)=\z(z)^n$.

   Recall that the polynomials   $P_{n}(z)$  are defined by equation  \e{eq:Jy} and the conditions  $P_{-1}(z)=0$,  $P_{0}(z)=1$.
   Put $ P(z)=\{ P_{n}(z)\}_{n=-1}^\infty$, $ f(z)=\{ f_{n}(z)\}_{n=-1}^\infty$,
    \begin{equation}
\omega (z):=  \{ P(z), f(z)\} = - 2^{-1}f_{-1}(z),
\label{eq:WR}\end{equation}
where the first formula \e{eq:Wr1}   has been used. In particular, for the operator $H_{0}$ we have  $ \omega_{0}(z) = -(2\z (z) )^{-1}$. By analogy with the continuous case, we define
     the Jost function
     \begin{equation}
 \Omega(z): = \omega(z) / \omega_{0}(z)= -2 \z(z) \omega (z).
\label{eq:RV}\end{equation}

The following result is a direct consequence of Theorem~\ref{Jost}.

\begin{corollary}\label{JOST}
Under  assumption   \e{eq:Tr} the Jost function $\Omega(z)$   depends analytically on $z\in\Pi$, and it is continuous in $z$ up to the cut along $[-1,1]$ except, possibly, the points $\pm 1$.  A point $z\in \Pi$ is an eigenvalue of the operator $H$ if and only if $\omega(z)=0$.
\end{corollary}

 Theorem~\ref{Jost} can be supplemented by the following assertion.
  
  \begin{theorem}\label{JostP}
  Under the assumptions of Theorem~\ref{Jost} put $a^{(N)}_{n}= a _{n}$, $b^{(N)}_{n}= b _{n}$ for $n \leq N$ and
  $a^{(N)}_{n}= 1/2$, $b^{(N)}_{n}=0$ for $n > N$. Let $  f_{n}^{(N)}(z)$ be the Jost solution  of equation \e{eq:Jy} with the coefficients
  $a^{(N)}_{n} $, $b^{(N)}_{n} $. Then for each $n\geq -1$, we have
 \begin{equation}
  \lim_{N\to\infty}  f_{n}^{(N)}(z)=f_{n} (z),\q z\in \clos\Pi , \q z\neq \pm 1.
\label{eq:cut}\end{equation}
\end{theorem}

\medskip

 {\bf 2.3.} 
 For $ \lambda\in (-1,1)$,    equation \e{eq:Jy} where $z=\lambda$ has two    solutions $f  (\lambda\pm i0)=\{f_{n} (\lambda\pm i0)\}_{n=-1}^\infty$;  of course $f_{n}(\lambda-i0)=\ov{f_{n} (\lambda+i0)}  $.  It follows from Theorem~\ref{Jost} and formula \e{eq:zt}  that
 asymptotics of  $  f_{n}(\lambda\pm i0)$ as $n\to\infty$  is oscillating:   
  \begin{equation}
  f_{n}(\lambda\pm i0)=e^{\mp i \theta n}  (1+O(\rho_{n})), \q \theta =\arccos\lambda\in (0,\pi).
\label{eq:Jost1}\end{equation}
Calculating their Wronskian \e{eq:Wr}  for $n\to\infty$, we find that
      \begin{equation}
\{f  (\lambda+i0),  f (\lambda-i0)\}= i\sin \theta\neq 0,\q   \sin\theta=\sqrt{1-\lambda^2},
\label{eq:HH1}\end{equation}
so that  these solutions  are linearly independent.  Thus,
 \begin{equation}
P_{n} (\lambda)= \ov{c(\lambda)} f_{n} (\lambda+i0)+ c(\lambda)  f_{n} (\lambda-i0)
\label{eq:HH2}\end{equation}
for some constant $c (\lambda)$. Taking the Wronskians of this equation with $f  (\lambda+i0)$   and using \e{eq:WR}, \e{eq:HH1}, we find that  
 \[
  \sin\theta \,c (\lambda) =i \omega(\lambda+i0)  .
\]
Of course $\omega(\lambda-i0)=\ov{\omega(\lambda+i0)}  $.
Thus \e{eq:HH2} leads to an intermediary result.

\begin{lemma}\label{HH}
  Under assumption \e{eq:Tr} the representation 
   \begin{equation}
P_{n} (\lambda)=\frac{ \omega(\lambda-i0)   f_{n} (\lambda+i0) -  \omega(\lambda+i0)  f_{n} (\lambda- i0)}{i\sqrt{1-\lambda^2}},  \q \lambda\in (-1,1), \q n=0,1,2, \ldots, 
\label{eq:HH4}\end{equation}
 holds true.
 \end{lemma}
 
In particular,   representation \e{eq:HH4} implies that
 \begin{equation}
\omega(\lambda\pm i0)  \neq 0, \q \lambda\in (-1,1).
\label{eq:HH5}\end{equation}
Indeed, otherwise we would have $P_{n} (\lambda)=0$ for some $\lambda\in (-1,1)$ and all $n\in{\Bbb Z}_{+}$.  However, 
 $P_0 (\lambda)=1$ for all $\lambda$.

 Let us set
    \begin{equation}
\varkappa (\theta) = | \Omega(\cos\theta+i0) |, \q \Omega(\cos\theta+i0)= \varkappa (\theta)e^{i \eta (\theta)},
\q \theta\in (0,\pi).
\label{eq:AP}\end{equation}
In the theory of the Schr\"odinger operator, the functions $\varkappa (\theta) $ and $\eta (\theta)$ are known as the limit amplitude and the limit phase, respectively; the function   $\eta (\theta)$ is also known as the scattering  phase or the   phase shift.  Definition \e{eq:AP}     fixes $\eta (\theta)$  only up to a term $2\pi k$ where $k\in{\Bbb Z}$.

It follows from  \e{eq:zt} and \e{eq:RV} that 
\[
2 \omega(\lambda\pm i0)=-\Omega(\lambda\pm i0) \z(\lambda\pm i0)^{-1}= - \varkappa (\theta) e^{\pm i(\eta (\theta)+\theta)}.
\]
Therefore, combined together  relations \e{eq:Jost1} and \e{eq:HH4} yield the Bernstein-Szeg\"o asymptotics of the  polynomials  $P_{n} (\lambda)$.
  
 \begin{theorem}\label{Sz}
  Under assumption \e{eq:Tr} for $\lambda \in (-1,1)$ the  polynomials $P_{n} (\lambda)$ have asymptotics 
   \begin{equation}
 P_{n} (\lambda)=  \varkappa (\theta) (sin\,\theta)^{-1} \sin ((n+1) \theta + \eta(\theta) ) + O(\rho_{n}) , \q \theta=\arccos \lambda,
\label{eq:Sz}\end{equation}
as $n\to\infty$. Relation \e{eq:Sz} is uniform in $\lambda$ on compact subintervals of $(-1,1)$.
 \end{theorem}

  \bigskip

 {\bf 2.4.}
 Let us construct
 the resolvent $R(z)=(H-zI)^{-1}$ of the operator $H$. Recall that   $\omega(z)$ is the Wronskian \e{eq:WR}.

 \begin{lemma}\label{res}
 For all $n,m\in{\Bbb Z}_{+}$, we have
  \begin{equation}
(R (z)e_{n}, e_{m})= \omega(z)^{-1} P_{n} (z) f_{m}(z),\q z\in \Pi, 
\label{eq:RR}\end{equation}
if $n\leq m$ and $(R (z)e_{n}, e_{m})=(R (z)e_m, e_n)$.
 \end{lemma} 
 
  \begin{pf} 
  We will show that
   the operator $R(z)$   defined by \e{eq:RR} is the resolvent of $H$. We have
    \begin{equation}
  \omega(z)( R (z)u)_{n} =  f_{n}(z) A_{n}(z)+   P_{n}(z) B_{n}(z)
\label{eq:RR1}\end{equation}
where
   \begin{equation}
A_{n}(z) =\sum_{m=0}^n  P_{m} (z) u_{m},\q    B_{n}(z) =\sum_{m=n+1}^\infty  f_{m} (z) u_{m},
\label{eq:RR2}\end{equation}
at least for all sequences $u=\{u_{n}\}$ with a  finite number of non-zero components $u_{n}$. In this case $R (z)u\in \ell^2 ({\Bbb Z}_{+})$ because $f_{n} (z)\in \ell^2 ({\Bbb Z}_{+})$ for all $z\in \Pi$. 

Our goal is to check that $(H-z)R (z)u=u$. 
It follows from definition \e{eq:ZP+}  of  the Jacobi operator $H$ and formula \e{eq:RR1} that
   \begin{multline}
  \omega ( (H-z) R  u)_{n} =  a_{n-1}\big(f_{n-1}A_{n-1}+P_{n-1} B_{n-1} \big)
  \\
  + (b_{n}-z) \big(f_{n}A_{n}+ P_{n} B_{n}  \big)+ a_{n}\big(f_{n+1}A_{n+1} +P_{n+1} B_{n+1}    \big).
\label{eq:RR3}
\end{multline}
According to \e{eq:RR2} we have
\[
f_{n-1}A_{n-1}+P_{n-1} B_{n-1}=f_{n-1}(A_{n}-P_{n} u_{n})+  P_{n-1}(B_{n}+f_{n} u_{n})
\]
and
\[
f_{n+1} A_{n+1} +P_{n+1} B_{n+1} =
f_{n+1}A_{n}  +  P_{n+1}B_{n} .
\]
Let us substitute these expressions into the right-hand side of \e{eq:RR3} and observe that
the coefficients at $A_{n}$ and $B_{n}$ equal zero by virtue of equation \e{eq:Jy} for $\{f_{n}\}$ and $\{P_{n}\}$, respectively. It  follows that
\[
 ( (H-z) R  u)_{n} = \omega^{-1}  a_{n-1}  (- P_{n}  f_{n-1} +    f_{n} P_{n-1}) u_{n}  =  u_{n} 
\]
whence  $R(z)=(H-z)^{-1}$; in particular, the operator $R(z) $ defined by \e{eq:RR} is bounded in $  \ell^2 ({\Bbb Z}_{+})$.
    \end{pf}

In view of Theorem~\ref{Jost},   $f_{n} (z)$ and $\omega (z)$ are continuous    functions of $z\in{\Bbb C}\setminus [-1,1]$ up to the cut along $[-1,1]$ with possible exception of the points $z=\pm 1$. Therefore using \e{eq:HH5}, we obtain the following result.

 \begin{theorem}\label{AC}
  Let assumption \e{eq:Tr} hold. Then, for all $n,m\in{\Bbb Z}_{+}$, the functions $(R(z)e_{n}, e_{m})$ are continuous as $z\in\Pi$ approaches the interval $(-1,1)$ from above or below, and the spectrum of the operator $H$ is absolutely continuous on $(-1,1)$.
 \end{theorem}
 
 We emphasize that the points $1$ and $-1$ may be eigenvalues of $H$; see Example~\ref{EE} below.
 
Let us now calculate  the spectral family  $d E(\lambda)$ of the   operator $H $. We proceed from the identity
\begin{equation}
2\pi i \frac{d(E (\lambda)e_n, e_m)} {d\lambda}=  (R (\lambda+ i0)e_n, e_m)-(R (\lambda- i0)e_n, e_m).
\label{eq:Priv}\end{equation}
It follows from formula \e{eq:RR} that
\[
(R (\lambda\pm i0)e_{n}, e_{m})=  \omega(\lambda\pm i 0)^{-1} P_{n} (\lambda) f_{m}(\lambda\pm i 0),\q n\leq m.
\]
Substituting this expression into \e{eq:Priv}, we find that
\[
2\pi i \frac{d(E (\lambda)e_n, e_m)} {d\lambda}= P_{n} (\lambda) \frac{ \omega(\lambda-i 0)  f_{m}(\lambda+i 0)    -\omega(\lambda+i 0) f_{m}(\lambda-i 0)    }{| \omega(\lambda+i 0) |^2}.
\]

Let us
combine this representation with  formula \e{eq:HH4} for $P_m (\lambda) $.
Since $ 2 | \omega(\lambda+i0) |=| \Omega(\lambda+i0) |$,
 we obtain the following result.

 \begin{theorem}\label{SF}
  Let assumption \e{eq:Tr} hold. Then, for all $n,m\in{\Bbb Z}_{+}$ and $\lambda\in(-1,1)$, we have the representation
  \[
 \frac{d(E (\lambda)e_n, e_m)} {d\lambda}= 2\pi^{-1}\sqrt{1-\lambda^2}  | \Omega(\lambda+i0) |^{-2} P_{n} (\lambda) P_m (\lambda)   .
\]
 \end{theorem} 
 
 \begin{corollary}\label{SFr}
 For $\lambda\in(-1,1)$,
the spectral measure of the operator $H$ equals
  \begin{equation}
d\rho(\lambda):=d(E (\lambda)e_0, e_0)= w(\lambda) d\lambda,
\label{eq:SFx}\end{equation}
where the weight function
 \begin{equation}
w(\lambda)=  2\pi^{-1}\sqrt{1-\lambda^2}  | \Omega( \lambda+i0) |^{-2}     .
\label{eq:SF1}\end{equation}
In particular, for the operator $H_{0}$, we have
\begin{equation}
w_{0} (\lambda) =2\pi^{-1}\sqrt{1-\lambda^2}.
\label{eq:ww}\end{equation}
\end{corollary}
 
In view of \e{eq:AP}, \e{eq:SF1} the amplitude in  \e{eq:Sz}   can  be written as
    \begin{equation}
 \varkappa (\theta)( \sin\theta)^{-1}= 2^{1/2} \pi^{-1/2}  (1-\lambda^2)^{-1/4} w(\lambda)^{-1/2}  
\label{eq:Sz1}\end{equation}
which is more common in the orthogonal polynomials literature.

Note that Theorem~\ref{Jost}  does not give any information on the behavior of the Jost function $\Omega(z)$ as $z\to\pm 1$. However relation \e{eq:SF1} implies that
 \begin{equation}
  \int_{-1}^1\sqrt{1-\lambda^2}  | \Omega( \lambda+i0) |^{-2} d\lambda
  = \frac{\pi}{2}
  \int_{-1}^1 w(\lambda) d\lambda\leq \pi/2,
\label{eq:WF}\end{equation}
and hence $ \Omega( \lambda+i0)$ cannot vanish too rapidly as $\lambda\to 1-0$ and $\lambda\to -1+0$ (even if $1$ or $-1$ are eigenvalues of $H$).
  For example, the behavior $\Omega( \lambda+i0)\sim c_{\pm}  (\lambda\mp 1)$ where $c_{\pm}  \neq 0$ is excluded.  

\section{Regular solutions of the Jacobi equations }

 {\bf 3.1.} 
 The regular solution $P(z)=\{P_{n}  (z)\}$ of the Jacobi equation \e{eq:Py} is determined by the boundary conditions $P_{-1}(z)=0$, $P_0(z)=1$. Then  $P_n(z) $ is a polynomial of degree $n$ and   relation \e{eq:J5} is satisfied. By analogy with the continuous case  (see, e.g., part~1 of Section~4.1  in  \cite{YA}), we here obtain bounds on  $P_n(z) $ for large $n$. As usual, the variable $\z(z)$ is defined by \e{eq:ome}.

 \begin{theorem}\label{OPP}
  If assumption \e{eq:Tr}  is satisfied, then
  \begin{equation}
|  P_{n}(z)|\leq C |\z (z)|^{-n} 
\label{eq:aa2}\end{equation}
with  some   positive constant $C$  not depending on $n$ and on $z$ in compact subsets of ${\Bbb C}\setminus\{-1,1\}$.
\end{theorem} 

Theorem~\ref{OPP} will be proven in Appendix~A.2.
 The proof relies on  the   equation (cf. equation \e{eq:V})
   \begin{equation}
 P_{n}(z)=P_{n}^{(0)}(z) +\frac{1}{\sqrt{ z^2-1}}\sum_{m=0}^{ n-1} (\z(z)^{n-m}- \z(z)^{m-n}) (VP(z))_{m},\q n\geq 1,
\label{eq:OPR1}\end{equation}
where $z\in\Pi$, $P_{n}^{(0)}(z)$ (normalized Chebyshev polynomials of the second kind) are given by equation \e{eq:OPR} and $V$ is defined by formula \e{eq:JV}.

Recall now  formula \e{eq:RR} for the resolvent $R(z)$ of the operator $H$.
Putting together estimates \e{eq:Jost} and \e{eq:aa2}, we obtain   the following result.

\begin{theorem}\label{ORP}
 Let  assumption   \e{eq:Tr} be satisfied. For all $n,m\in  {\Bbb Z}_{+}$, the functions $(R (z)e_{n}, e_{m})$ 
 are  continuous in $z$ up to the cut along $[-1,1]$ except, possibly, the points $\pm 1$.
 Moreover,
  \[
|(R (z)e_{n}, e_{m})|\leq C|\omega(z)|^{-1}   |\z (z)|^{|n-m|}  
\]
with  some   positive constant $C $ that does not depend on $n$, $m$ and  on $z$ in compact subsets of $\clos\Pi $ and away from the points $\pm 1$.
\end{theorem}

We need also a representation of the Jost function $\Omega(z)$ in terms of the orthogonal polynomials $P_{n}  (z)$. It plays the role of the representation (see, e.g., formula (1.38) in Chapter~4 of \cite{YA}) of the Jost function via the regular solution of the Schr\"odinger equation in the continuous case.

\begin{proposition}\label{ORPj}
 Let  assumption   \e{eq:Tr} be satisfied, and let the Jost function $\Omega(z)$ be defined by \e{eq:WR}. Then
   \begin{equation}
\Omega(z)=1-2\sum_{n=0}^{\infty} \z(z)^{n+1} (VP(z))_{n} 
\label{eq:OM}\end{equation}
for all $z\in\clos\Pi$  except, possibly, the points $\pm 1$.
\end{proposition}
 
 \begin{pf}
 Suppose first that $z=\lambda + i0$ where $\lambda\in (-1,1)$.
Substituting expressions \e{eq:OPs} and \e{eq:HH4}   into  equation \e{eq:OPR1} we see  that
 \begin{multline*}
 \omega(\lambda-i0)   f_{n} (\lambda+i0) -  \omega(\lambda+i0)  f_{n} (\lambda- i0)=2^{-1} (e^{i(n+1)\theta}-
 e^{-i(n+1)\theta})
 \\
 +\sum_{m=0}^{ n-1} (e^{-i(n-m)\theta}- e^{i(n-m)\theta}) (VP(\lambda + i0))_{m}.
\end{multline*}
Let us consider the asymptotics of both sides of this equation as $n\to\infty$ and use relation \e{eq:Jost1}. Then comparing coefficients at $e^{i n\theta} $, we get \e{eq:OM} for 
$z=\lambda + i0$. Both sides of \e{eq:OM} are analytic in $z\in \Pi$ and continuous up to the cut along $[-1,1]$ according to Theorems~\ref{Jost} and \ref{OPP}.  Therefore relation \e{eq:OM} extends to all $z$.
 \end{pf}

\bigskip

 {\bf 3.2.}  
 Let us   find    asymptotics  of the polynomials $P_{n} (z)$  for $z\not\in [-1,1]$. We follow here closely the scheme exposed in Section~4.1 (see, in particular, Lemma~1.11) of \cite{YA}. 
 
 We start by introducing solutions $g_{n} (z)$ of equation  \e{eq:Jy} exponentially growing as $n\to\infty$. Perhaps this construction is of interest in its own sake. For $z\in \Pi$,  fix $n_{0}= n_{0} (z)$ such that $f_{n} (z)\neq 0$ for $n\geq n_{0}-1$. Note that, for $\Im z\neq 0$,    one can set $n_{0}= 0$ because the equality $f_{n_{0}-1} (z)= 0$  implies that the Jacobi operator $H^{(n_0)}$ with the matrix elements $a^{(n_0)}_{n}=  a_{n+n_{0}}$, $b^{(n_0)}_{n}=  b_{n+n_{0}}$ has the eigenvalue $z$.
  Put
   \begin{equation}
g_{n} (z)= f_{n}(z)\Theta_{n} (z)\q {\rm where} \q \Theta_{n} (z)=\sum_{m=n_{0}}^n (a_{m-1} f_{m-1}(z) f_{m}(z))^{-1}.
\label{eq:GE}\end{equation}
An elementary calculation shows that this sequence satisfies equation \e{eq:Jy}. It is also easy to find the asymptotics of $g_{n}(z)$ as $n\to\infty$:
 \begin{multline*}
\z^n g_{n}(z)=2\z^{2n}\sum_{m=n_{0}}^n \z ^{-2m+1}(1+o(1)) = 2\z^{2n+1} \frac{\z^{-2n-2}-1}{\z^{-2}-1}(1+o(1))
\\
= 2 \frac{ 1-\z^{2n+2}}{\z^{-1}-\z}(1+o(1)) =  \frac{1+o(1)}{\sqrt{z^2-1}} 
\end{multline*}
where $\z=\z(z)$.
This yields the following result.

\begin{lemma}\label{GE}
Let   $z\in \Pi$. 
  Under assumption \e{eq:Tr} the sequence $g_{n} (z)$ defined by \e{eq:GE} satisfies equation \e{eq:Jy} and
   \[
\lim_{n\to\infty} \z(z)^n g_{n}(z)=\frac{1}{\sqrt{z^2-1}} .
\]
 \end{lemma}
 
 The Wronskian \e{eq:Wr}  of $f (z)=\{f_{n}(z)\}$, $g (z)=\{g_{n}(z)\}$ equals
  \[
\{ f(z),g(z)\}= a_{n}f_{n}(z)f_{n+1}(z)( \Theta_{n+1}(z)- \Theta_{n}(z))=1,
\]
and hence solutions $f(z)$ and $g (z)$ are linearly independent. It follows that 
 \begin{equation}
P_{n} (z)= d_{+}(z)f_{n} (z)+d_{-}(z)g_{n} (z)
\label{eq:GE4}\end{equation}
where 
\[
\{ P(z),f(z)\}=  d_{-}(z)\{ g(z), f(z)\}= -d_{-}(z)
\]
 and 
 \[
 \{ P(z),g(z)\}=  d_{+}(z)\{ f(z), g(z)\}= d_{+}(z). 
 \]
 According to  \e{eq:WR}, \e{eq:RV} we have 
  \begin{equation}
d_{-}(z)=(2\z (z))^{-1}\Omega  (z).
\label{eq:cba}\end{equation}
  Obviously,  
$d_{+}(z)\neq 0$ if $d_{-}(z)= 0$. Therefore Lemma~\ref{GE} implies the following result.

\begin{theorem}\label{GE1}
 Under assumption \e{eq:Tr} for all  $z\in \Pi$, we have the relation 
   \begin{equation}
\lim_{n\to\infty} \z(z)^n P_{n}(z)=\frac{\Omega(z)}{1-\z(z)^2}
\label{eq:GEGE}\end{equation}
with convergence uniform on compact subsets of $\Pi$.
Moreover, if $\Omega(z)=0$, then
 \begin{equation}
\lim_{n\to\infty} \z(z)^{-n} P_{n}(z)=\{ P(z),g(z)\} \neq 0.
\label{eq:GEGEx}\end{equation}
 \end{theorem}
 
The existence of the limit in \e{eq:GEGE} is the  classical result of the Szeg\"o theory. It is stated as Theorem~12.1.2 in the book  \cite{Sz} where the assumptions are imposed   on the measure $d\rho(\lambda)$; in particular, it is assumed that
$\supp\rho\subset [-1,1]$. Under assumption \e{eq:Tr} asymptotic relation \e{eq:GEGE} was established in Theorem~2.20 of \cite{KS}; our proof is rather different from that in \cite{KS}. Relation \e{eq:GEGEx} is perhaps new.

 \section{ Edge points of the spectrum}  
 
   {\bf 4.1.} 
   A study of  the Jost function $\Omega(z)$ as  $z\to 1$ and $z\to -1$ requires an additional assumption on the coefficients $a_{n}$, $b_{n}$. This is quite similar to the Schr\"odinger operator in the space $L^{2}({\Bbb R}_{+})$ that has a threshold at the point $z=0$. 
 The presentation here follows   very closely Section~4.3 of \cite{YA} where the continuous case was considered; so some technical details will be omitted.

\begin{theorem}\label{JOSF}
Suppose that 
 \begin{equation}
\sum_{n=0}^\infty n(|  a_{n}-1/2|+|b_{n}|)<\infty. 
\label{eq:Trfd}\end{equation} 
Then all functions $f_{n} (z)$ and, in particular, 
 the Jost function $\Omega(z)$  are continuous in $z$ up to the cut along $[-1,1]$. Moreover, $f_{n}(\pm 1)$ satisfies  equation \e{eq:Jy} for $z=\pm 1$, that is
  \begin{equation}
 a_{n-1} f_{n-1} (\pm 1) +b_{n} f_{n} (\pm 1) + a_{n} f_{n+1} (\pm 1)= \pm f_{n} (\pm 1), \q n\in{\Bbb Z}_{+},
\label{eq:Jypm}\end{equation}
 and  $f_{n}(\pm 1)=(\pm 1)^n +o(1)$ as $n\to\infty$.
\end{theorem}

Theorem~\ref{OPP} can be supplemented by the following result.

 \begin{theorem}\label{OPP1}
  If assumption \e{eq:Trfd}  is satisfied, then
  \begin{equation}
|  P_{n}(z)|\leq C (n+1) |\z (z)|^{-n} 
\label{eq:aa2+}\end{equation}
with  some   positive constant $C$  not depending on $n$ and on $z$ in compact subsets of ${\Bbb C} $.
\end{theorem}

The proofs of Theorem~\ref{JOSF} and \ref{OPP1} will be discussed in Appendix~A.3. We emphasize that assumption \e{eq:Trfd} cannot be relaxed; see Remark~\ref{find1} below.

Passing in \e{eq:V} to the limit $z\to\pm 1$ and taking into account the relation
\begin{equation}
\z(z)^{k} =(\pm 1)^{k} (1 \mp k \sqrt{z^{2}-1}) + O ( |z^{2}-1|), \q z\to\pm 1,
\label{eq:jj}\end{equation}
 we obtain an equation for the sequences $f_{n}(\pm 1)$:
  \[
  f_{n}(\pm 1)=(\pm 1)^n +2 \sum_{m=n+1}^\infty (\pm 1)^{n-m+1}( n-m)  (Vf(\pm 1))_{m}.
\]
Similarly, since $P_{n}^{(0)}(\pm 1)= (n+1) (\pm 1)^{n}$, equation \e{eq:OPR1}  yields
 \begin{equation}
  P_{n}(\pm 1)= (n+1) (\pm 1)^n -2 \sum_{m=0}^{n-1} (\pm 1)^{n-m+1}( n-m)  (V P(\pm 1))_{m}.
\label{eq:Vpmr}\end{equation}

Estimate \e{eq:aa2+} allows us to
  pass  to the limit $z\to\pm 1$ in \e{eq:OM} which leads to the following result.

\begin{lemma}\label{ORPj1}
 Under   assumption   \e{eq:Trfd}, the representation
   \begin{equation}
\Omega(\pm 1)=1-2\sum_{n=0}^{\infty} (\pm 1)^{n+1} (VP(\pm 1))_{n} 
\label{eq:OM1}\end{equation}
holds.
\end{lemma}

Similarly to  Lemma~\ref{GE}, we can introduce a solution  $g_{n} (\pm 1)$ of the equation
\e{eq:Jypm} linearly independent with  $f_{n} (\pm 1)$.

 \begin{lemma}\label{GEpm}
Let  assumption \e{eq:Trfd} hold. Define the sequence $g_{n} (\pm 1)$ by formula \e{eq:GE}. Then $g_{n} (\pm 1)$ satisfies equation \e{eq:Jypm}, $g_{n} (\pm 1)= 2(\pm 1)^{n +1}n (1+o(1))$ as $n\to\infty$ and $\{f(\pm 1),
g(\pm 1)\}=1$. 
\end{lemma}

Since neither of solutions  $f_{n} (\pm 1)$,  $g_{n} (\pm 1)$ nor their linear combinations   tend to zero as $n\to \infty$, we obtain

\begin{theorem}\label{JOSF2}
Under  assumption \e{eq:Trfd} equations \e{eq:Jypm} do not have solutions tending to zero as $n\to \infty$. In particular, the operator $H$ cannot have eigenvalues $1$ and $-1$.
\end{theorem}

To find the asymptotics of the polynomials $P_{n} (z)$ at critical points $z=\pm 1$, we use equalities \e{eq:GE4} and \e{eq:cba} where $z=\pm 1$. The next result is a  direct consequence of Lemma~\ref{GEpm}.

\begin{theorem}\label{AScr}
Let  assumption \e{eq:Trfd}, and let  $\Omega (\pm 1) \neq 0$. Then
\begin{equation}
P_{n} (\pm 1)= \Omega (\pm 1) (\pm 1)^{n  }n  (1+o(1))
\label{eq:P+}\end{equation}
\end{theorem}

 Of course, if $\Omega (\pm 1) \neq 0$, then according to  \e{eq:RR} the resolvent kernel $(R(z) e_{n}, e_{m})$ is continuous as $z\to \pm 1$.

\bigskip

   {\bf 4.2.} 
    Our next goal is to find the asymptotic behavior as $z\to\pm 1$ of $\Omega (z)  $ and hence of $(R(z) e_{n}, e_{m})$
in the exceptional case $\Omega (\pm 1) = 0$. Let us use the same terminology as for Schr\"odinger operators.    

   \begin{definition}\label{resin}
Let  assumption \e{eq:Trfd} hold. If $\Omega (\pm 1) = 0$, we say that the operator $H$ has a resonance at $z=\pm 1$.
\end{definition}

Clearly, the condition  $\Omega (\pm 1) = 0$ is equivalent to the linear dependence of the solutions $P_{n} (\pm 1)$ and
$f_{n} (\pm 1)$ of equation \e{eq:Jypm}. 
In this case the sequence $P_{n}(\pm 1)$ is bounded as  $n\to\infty$. The next result follows from equations
\e{eq:Vpmr} and \e{eq:OM1}.

 \begin{lemma}\label{resin2}
  Suppose that $\Omega (\pm 1) = 0$ and put
 \begin{equation}
  \gamma_{\pm}= 1 + 2\sum_{m=0}^{\infty} (\pm 1)^{m-1}  m(VP (\pm 1))_{m}.
\label{eq:Vx2}\end{equation}
Then   there exists
 \begin{equation}
  \lim_{n\to\infty} (\pm 1)^{n}  P_{n}(\pm 1)=\gamma_{\pm} \neq 0.
\label{eq:Vx1}\end{equation}
In particular, $  P_{n}(\pm 1)=  (\pm 1)^{n}  \gamma_{\pm} f_{n}(\pm 1)$.
\end{lemma}

Now we  are in a position to find the asymptotic behavior of the Jost function $\Omega (z)  $ as $z\to \pm 1$ in the case $\Omega (\pm 1) = 0$. Let us proceed from representation \e{eq:OM}. According to  \e{eq:jj}, for each $n$,   we have
\[
\z(z)^{n+1}(VP(z))_{n} =  \big((\pm 1)^{n +1}- (\pm 1)^{n } (n+1)\sqrt{z^{2}-1} \big)(VP(\pm 1 ))_{n}  + O (z^2-1).
\]
Substituting this expression into \e{eq:OM} and taking into account equation \e{eq:OM1},  we see that
 \[
\Omega(z)=  \Omega(\pm 1) + 2 \sum_{n=0}^\infty (\pm 1)^{n } (n+1) (VP(\pm 1 ))_{n}
  \sqrt{z^{2}-1} + o (\sqrt{|z^{2}-1|}).
\]
If $\Omega (\pm 1) = 0$, then the coefficient at $  \sqrt{z^{2}-1}$ here equals $ \pm \gamma_{\pm}$. Let us state the result obtained.

\begin{theorem}\label{resin1}
Under  assumption \e{eq:Trfd} suppose that $\Omega (\pm 1) = 0$ and define the number $  \gamma_{\pm}$  by \e{eq:Vx2} or, equivalently, by \e{eq:Vx1}. Then
 \[
\Omega(z)=\pm \gamma_{\pm}  \sqrt{z^{2}-1} + o (\sqrt{|z^{2}-1|})
\]
as $z\to \pm 1$.
\end{theorem}

Using formula   \e{eq:SF1}, we obtain the following consequence for the weight function.

\begin{corollary}\label{resnw}
Under the  assumptions of Theorem~\ref{resin1}, we have
 \[
w(\lambda)=2\pi^{-1}\gamma_{\pm}^{-2}(1-\lambda^2)^{-1/2} (1+ o(1)), \q \lambda\in (-1,1),
\]
as $\lambda\to \pm 1$.
\end{corollary}

 In view of Lemmas~\ref{res} and \ref{resin2}, Theorem~\ref{resin1} implies also the following result.

\begin{corollary}\label{resn1}
Under the  assumptions of Theorem~\ref{resin1}  for all $n,m\in{\Bbb Z}_{+}$, the representation 
 \[
 (R(z)e_{n}, e_{m})= -2\frac{ (\pm 1)^{\min\{n,m\} }f_{n} (\pm 1) f_{m} (\pm 1) + o(1)}
{\sqrt{z^{2}-1}  }
\]
as $z\to \pm 1$ is satisfied.
\end{corollary}

\begin{remark}\label{find}
Combined together, Theorems~\ref{JOSF} and \ref{resin1} ensure that
under  assumption \e{eq:Trfd} the discrete spectrum of the operator $H$ is finite.
\end{remark}

\begin{remark}\label{find1}
Under assumption   \e{eq:Trfd} the asymptotics of the orthogonal polynomials $P_{n} (\pm 1)$ is given by formulas
 \e{eq:P+} and  \e{eq:Vx1}. An example  of the Jacobi polynomials shows   that the asymptotics of  $P_{n} (\pm 1)$ is significantly   changed  (see  formula \e{eq:th}, below) if the condition \e{eq:Trfd} is even slightly relaxed. In particular, estimate \e{eq:aa2+} is violated in this case.
 \end{remark}

\begin{remark}\label{find2}
The problem of characterization of spectral data corresponding to assumption   \e{eq:Trfd} was considered in \cite{Gus}.
 \end{remark}
  
   \bigskip

 {\bf 4.3.}
  The condition \e{eq:Trfd} of Theorem~\ref{JOSF2} guaranteeing, in particular,  the absence of eigenvalues $\pm 1$ is optimal. Indeed, 
    there exist Jacobi operators satisfying condition \e{eq:Tr} but not \e{eq:Trfd}  with an eigenvalue at the points $+1$ or $-1$. In  the example below, we suppose that $a_{n}=1/2$.  Set
   \begin{equation}
 \psi_{n}^{(\pm)}=(\pm 1)^n (n+1)^{-l} 
\label{eq:EE1}\end{equation}
for all $n\in {\Bbb Z}_{+}$ and $ \psi_{-1}^{(\pm)}=0$. 
Then $\psi^{(\pm)}=\{\psi_{n}^{(\pm)}\}\in \ell^2 ({\Bbb Z}_{+})$ if $l>1/2$.  Let us consider  equation   \e{eq:Jy} where
$u_{n}=\psi_{n}^{(\pm)}$ and $z=\pm 1$   as
 an equation for $b_{n}^{(\pm)}$:
     \begin{equation}
 b_{n}^{(\pm)}= \pm 1- (2 \psi_{n}^{(\pm)})^{-1}( \psi_{n-1}^{(\pm)}+ \psi_{n+1}^{(\pm)}) .
\label{eq:EE2}\end{equation}
It follows from \e{eq:EE1} that   
  \[
  b_{n}^{(\pm)}  = \mp 2^{-1}l(l+3) n^{-2}+O (n^{-3})  
\]
as $n\to\infty$. In particular, condition \e{eq:Tr} holds. 
 Moreover, up to a finite rank operator, the perturbation $V^{(\pm)}=H^{(\pm)}-H_{0}$ of the operator $H_{0}$ is negative for the upper sign   and  positive for the lower sign. 
 Thus, we have
 
 \begin{example}\label{EE}
 Let  $a_{n}=1/2$ and $b_{n}^{(\pm)}$ be given by formulas \e{eq:EE1}  where $ l>1/2 $ and \e{eq:EE2}.
Then condition \e{eq:Tr}  is satisfied. The Jacobi operator $H^{(\pm)}$   has the eigenvalue $\pm 1$ with the eigenvector $\psi_{n}^{(\pm)}$ defined by \e{eq:EE1}.  The spectrum of $H^{(+)} $ is finite above the point $1$, and spectrum of $H^{(-)} $ is finite  below the point $-1$.
\end{example}

Note that   the operators $H^{(\pm)}$ have infinite number of eigenvalues. For instance, this fact follows from a general result of \cite{D-K}, Theorem~3, stating that the spectrum   of a Jacobi operator is purely absolutely continuous on $[-1,1]$ if it has only a finite number of eigenvalues outside $[-1,1]$. In our case $H^{(\pm)}$  has the eigenvalue $\pm 1$.

Example~\ref{EE}   exhibits an asymptotic behavior of orthogonal polynomials $P_{n} (z) $ intermediary between exponential  (for $z\not \in [-1,1]$) and oscillating (for  $z \in (-1,1)$) regimes. Indeed,  if $a_{n}=1/2$ and $b_{n}^{(\pm)}$ are defined by \e{eq:EE1} and  \e{eq:EE2}, then $P_{n}(\pm 1)=\psi_{n}^{(\pm)}$ are given by \e{eq:EE1}. So, this sequence behaves as some power  of $n$ (negative or positive)  as $n\to\infty$.

\bigskip
 
   {\bf 4.4.} 
A similar result is  true  for the Schr\"odinger operator  ${\bf H}=-d^2/dx^2+ b(x) $ with  the boundary condition $u(0)=0$
in the space $L^2 ({\Bbb R}_{+})$. Let $\psi\in C^\infty ({\Bbb R}_{+})$, $\psi(x)>0$ for all 
 $x>0$, $\psi(x)=x$ in a neighborhood of the point $x=0$ and
  \begin{equation}
 \psi(x)=x^{-l}\q {\rm where}\q l>1/2
\label{eq:EE1C}\end{equation}
for sufficiently  large $x$. Then $\psi \in L^2 ({\Bbb R}_{+})$.
Put
  \begin{equation}
 b(x)= \psi''(x) \psi(x)^{-1}
\label{eq:EE2C}\end{equation}
so that
 $
 b(x)=l(l+1) x^{-2}
$
if $x$ is large enough. This yields

 \begin{example}\label{EEC}
 Let  $b(x)$   be given by formulas \e{eq:EE1C} and \e{eq:EE2C}.
Then the  operator ${\bf H}$   has the eigenvalue $0$ with the eigenfunction  $\psi (x)$ defined by \e{eq:EE1C}, and its negative spectrum   is finite.
\end{example}

The decay of $b(x)$ as $|x|^{-2}$ at infinity is critical. Actually, it is   known (see part~1 in Section~4.3  of \cite{YA}) that 
${\bf H}$ cannot have the zero eigenvalue provided 
$$\int_{0}^\infty (1+x )|b(x)| dx<\infty.
$$ 

 \section{ The perturbation determinant and the spectral shift function. Trace identities}

 {\bf 5.1.} 
 Let $H_{0}$ and $H$ be arbitrary self-adjoint operators with a trace class difference  $V =H-H_{0}$. Then
 the perturbation  determinant
   \begin{equation}
\Delta (z)=\det \big(I+V R_{0} (z)\big)
\label{eq:Tr1}\end{equation} 
for the pair $H_{0}$, $H$
 is well defined and is an analytic function of $z\in{\Bbb C}\setminus\sigma (H_{0})$. 
 Obviously, $\Delta (\bar{z})=\overline{\Delta (z)}$ and
  \begin{equation}
 \Delta (z)\to 1\q {\rm as}\q \dist\{z,\sigma (H_{0}) \}\to\infty.
\label{eq:PDhe}\end{equation}
Note also  the general formula
  \begin{equation}
\tr \big( R(z)-R_{0}(z)\big)=- \frac{   \Delta' (z)  }{  \Delta (z)   }.
\label{eq:PD}\end{equation}
In view of our applications to Jacobi operators, below we suppose that the operators $H_{0}$ and $H$ are bounded.

  The  Kre\u{\i}n spectral shift function $\xi (\lambda)$ 
   is defined in terms of the perturbation determinant  \e{eq:Tr1}. According to \e{eq:PDhe} we can fix the branch of the function $\ln\Delta (z) $ for $\Im z\neq 0$ by the condition
  \begin{equation}
 \arg\Delta (z)\to  0 \q {\rm as}\q \dist\{z,\sigma (H_{0}) \}\to\infty.
\label{eq:PHE}\end{equation}
Then
   \begin{equation}
  \xi (\lambda):=\pi^{-1}\lim_{\varepsilon\to+0}\arg \Delta(\lambda+i\varepsilon).
\label{eq:SSF}\end{equation}
This limit exists for a.e. $ \lambda\in{\Bbb R}$,
\begin{equation}
 \int_{-\infty}^\infty |\xi(\lambda) |d\lambda \leq \| V\|_{1}
\label{eq:DS}\end{equation}
and the representation  
\begin{equation}
\ln\Delta (z) =\int_{-\infty}^\infty \xi(\lambda) (\lambda-z)^{-1} d\lambda, \q \Im z\neq 0,
\label{eq:DD}\end{equation}  
holds.  The    function $\xi (\lambda)$ is constant on intervals not containing points of $  \sigma(H_{0})\cup \sigma(H) $. In particular,  $\xi(\lambda)=0$ for $\lambda<\inf\big( \sigma(H_{0})\cup \sigma(H)\big)$ and $\lambda>\sup \big( \sigma(H_{0})\cup \sigma(H)\big)$. If $\lambda_{1}$ is an isolated eigenvalue of finite multiplicity $k_{0}$  of the operator $H_{0}$ and  multiplicity $k $  of the operator $H$, then 
\begin{equation}
\xi (\lambda_{1} +0)- \xi (\lambda_{1} -0)=k_{0}-k.
\label{eq:ES}\end{equation}

 We refer  to the books \cite{GK, Sim,Ya}   for a detailed presentation of  these notions.
 
 \bigskip
  
  {\bf 5.2.} 
  Let us come back to Jacobi operators.
 Under the assumption \e{eq:Tr}, the    perturbation  $V$    belongs to the trace class $\mathfrak{S}_{1}$ so that all results mentioned in the previous subsection are true.
 
   Let us find a link  of the perturbation determinant $\Delta(z)$ with the Jost function \e{eq:RV}. We follow here closely the scheme of the paper \cite{BF}  where the   operator $-D^2+ b(x)$ was considered;  for its detailed presentation, see   Section~4.1 of \cite{YA}. Note that there exists also a different       to proofs of such assertions relying on the Fredholm expansion of determinants; this approach was developed in  \cite{J-P}.
   
   Let us find an expression of the left-hand side of \e{eq:PD} in terms of the Jost function $\Omega(z)$. 
The first assertion is true without any assumptions on $a_{n}$ and $b_{n}$. Recall that $\{f_{n} (z)\}_{n=-1}^{\infty}$
is the Jost solution of equations \e{eq:Jy} and $\omega(z)$  is defined by \e{eq:WR}.

\begin{lemma}\label{RR1}
For all $z\in \Pi$ and all $N\geq 0$, we have the identity 
  \begin{equation}
\sum_{n=0}^N  f_{n}(z)P_{n } (z)  = - \omega'(z)  + a_{N} \big( P_{N}  (z) f_{N+1}'  (z)
- P_{N+1}  (z) f_{N }'  (z)\big). 
\label{eq:PP1}\end{equation}
 \end{lemma}

 \begin{pf}
 Let us differentiate (in $z$) equation \e{eq:Jy} for $f_{n} (z)$   and then multiply it by $P_{n} (z)$:
  \begin{equation}
 a_{n-1}f_{n-1}'(z)P_{n } (z)  + a_{n} f_{n+1}' (z) P_{n } (z) = (z-b_{n}) f_{n}' (z)P_{n } (z)+ f_{n} (z) P_{n } (z) .
\label{eq:Jy1}\end{equation}
Multilplying equation \e{eq:Py} for $P_{n} (z)$  by $f_{n}' (z)$, we see that
 \begin{equation}
 a_{n-1}P_{n-1}(z)f_{n }' (z)  + a_{n} P_{n+1} (z) f_{n }' (z) = (z-b_{n}) P_{n } (z) f_{n}' (z) .
\label{eq:Jy2}\end{equation}
Then we subtract \e{eq:Jy2} from \e{eq:Jy1}:
 \begin{multline}
f_{n} (z) P_{n } (z) =  \big(a_{n-1}f_{n-1}'(z)P_{n } (z)-a_n f_{n}'(z)P_{n+1 } (z)\big)
\\
+
 \big(a_{n}f_{n+1}'(z)P_{n } (z) -a_{n-1}f_{n}'(z)P_{n -1} (z) \big).
\label{eq:Jy3} \end{multline}

Let us calculate the sums of the right-hand sides over $n=0,1, \ldots, N$:
 \[
\sum_{n=0}^N \big(a_{n-1}f_{n-1}'(z)P_{n } (z)-a_n f_{n}'(z)P_{n+1 } (z)\big) = 
2^{-1} f_{-1}'(z)P_0 (z)-a_N f_N '(z)P_{N+1 } (z)  
\]
and 
 \[
\sum_{n=0}^N \big(a_{n}f_{n+1}'(z)P_{n } (z) -a_{n-1}f_{n}'(z)P_{n -1} (z)  \big) = a_{N}
  f_{N+1} '(z)P_{N  } (z) - 2^{-1} f_{0} '(z)P_{-1  } (z)   .
\]
Since $P_{-1} (z)=0$,   $P_{0} (z)=1$ and $f_{-1}(z)=-2 \omega(z)$, taking the sum of equations \e{eq:Jy3}, we obtain 
the identity \e{eq:PP1}.
   \end{pf}

  Let us calculate the asymptotics of the right-hand side of  \e{eq:PP1} as $n\to\infty$ supposing first that the perturbation $V$ has finite support, that is,
   $a_{n}=1/2$   and  $b_{n}= 0$ for sufficiently large $n$.  In view of   \e{eq:omex},  the equations \e{eq:Jy} are satisfied    in this case for all $ z\in \Pi$ and $n$  large enough
 if $u_n (z) $ is a  arbitrary linear combination of the functions $ \z (z)^n$ and $ \z (z)^{-n}$. In particular, we have
  \begin{equation}
f_{n} (z)=\z (z)^n \q{\rm and}\q P_{n}  (z)= \gamma_{+} (z) \z (z)^n+ \gamma_{-} (z) \z (z)^{-n}
\label{eq:PP}\end{equation}
for some numbers $\gamma_\pm (z)$. Recall that $|\z(z)|<1$ for $z\in\Pi$.
   
   \begin{lemma}\label{RR2}
For all $z\in \Pi$,  $z\not\in \sigma_{d} (H)$ and $n\to\infty$, we have the relation 
  \begin{equation}
  \omega(z)^{-1}\big( P_{n}  (z) f_{n+1}'  (z)
- P_{n+1}  (z) f_{n }'  (z)\big) = -\frac{ 2n  }{\sqrt{z^2-1}}+  \frac{ \z (z)}{ z^2-1}+O (\z(z)^{2n}).
\label{eq:PP2}\end{equation}
 \end{lemma}

 \begin{pf}
 It follows from formulas \e{eq:PP} that the Wronskian \e{eq:WR} calculated for large $n$ is given by the equation
  \begin{equation}
2 \omega=     (\gamma_{+}\z^n + \gamma_{-}\z^{-n}) \z^{n+1} -(\gamma_{+}\z^{n+1} + \gamma_{-}\z^{-n-1}) \z^{n} =
  \gamma_{-} (\z-\z^{-1})  .
\label{eq:PP3}\end{equation}
 Formulas \e{eq:PP}  also imply that
 \begin{multline}
 P_{n}f_{n+1}' - P_{n+1}f_{n }' 
=  (\gamma_{+}\z^n + \gamma_{-}\z^{-n}) (\z^{n+1})' - (\gamma_{+}\z^{n+1} + \gamma_{-}\z^{-n-1}) (\z^n)'
\\
= -\frac{1 }{\sqrt{z^2-1}}\Big( \gamma_{-}\big((n+1) \z- n  \z^{-1}\big)+\gamma_{+}\z^{2n+1}\Big)   
\label{eq:PP4}\end{multline}
where we used that
\[
\z (z)'=-\frac{\z (z)}{\sqrt{z^2-1}} .
\]
Dividing now \e{eq:PP4} by \e{eq:PP3}, we arrive at  \e{eq:PP2}.
  \end{pf}

Using expression \e{eq:RR} for $ (R (z)e_{n}, e_{n})$ and putting together Lemmas~\ref{RR1} and \ref{RR2},     we get the following result.
  
  \begin{theorem}\label{RRR}
  Suppose  that $a_{n}=1/2$   and  $b_{n}= 0$ for sufficiently large $n$. Then
for all $z\in \Pi$,  $z\not\in \sigma_{d} (H)$  and $N\to\infty$, we have the relation 
  \begin{equation}
\sum_{n=0}^N (R (z)e_{n}, e_{n})=- \omega'(z) \omega(z)^{-1}- \frac{ N  }{\sqrt{z^2-1}}+  \frac{ \z (z) }{2( z^2-1)}+O (\z(z)^{2N})  .
\label{eq:PPP}\end{equation}
 \end{theorem}
 
 Of course the same formula is true for the operator $H_{0}$ when $\omega_{0}(z)= -2^{-1} \z (z)^{-1}$.
 Comparing formulas \e{eq:PPP} for $H$ and $H_{0}$,  we find that
   \[
\sum_{n=0}^N\big( (R (z)e_{n}, e_{n})-(R_{0} (z)e_{n}, e_{n})\big)=- \omega'(z) \omega(z)^{-1} +  \omega_{0}'(z) \omega_{0}(z)^{-1} +O (\z(z)^{2N})  .
\]
Passing here to the limit $N\to\infty$, we see that 
  \begin{equation}
\tr \big( R(z)-R_{0}(z)\big)=- \frac{ \Omega'(z)} {   \Omega(z)},\q z \in \Pi,
\label{eq:PPP2}\end{equation}
where the function $\Omega(z)$  is defined by \e{eq:RV}.

Theorem~\ref{JostP} allows us to extend this result to arbitrary operators satisfying condition \e{eq:Tr}. Indeed,
set $a^{(N)}_{n}= a _{n}$, $b^{(N)}_{n}= b _{n}$ for $n \leq N$ and
  $a^{(N)}_{n}= 1/2$, $b^{(N)}_{n}=0$ for $n > N$. Let
 $H^{(N)}$ be the Jacobi operator with the  matrix elements $a^{(N)}_{n}$, $b^{(N)}_{n}$,  and let $\Omega^{(N)}(z)= - 2\z(z)\{P^{(N)}(z), f^{(N)}(z)\}$. Write formula \e{eq:PPP2} for the pair $H_{0}, H^{(N)}$ and pass to the limit $N\to\infty$. Since $V^{(N)}:=H^{(N)}-H_{0}\to V$ in the trace norm, we  see that
\[
\tr \big( R^{(N)}(z)-R_{0}(z)\big)\to  \tr \big( R(z)-R_{0}(z)\big)
\]
as $N\to \infty$. According to Theorem~\ref{JostP} we also have $\Omega^{(N)}(z)\to \Omega(z)$ and hence 
$d\Omega^{(N)}(z)/dz\to d\Omega(z)/dz$ as $N\to \infty$. This leads to the desired result.

 \begin{theorem}\label{RRX}
The relation \e{eq:PPP2} holds true under assumption \e{eq:Tr}.
 \end{theorem}

   \bigskip

   {\bf 5.3.} 
   Putting together relations \e{eq:PD} and  \e{eq:PPP2}, we see that
   \begin{equation}
\Delta (z)=   A \, \Omega  (z), \q z\in \Pi ,
\label{eq:HH}\end{equation}
for some constant $A\in{\Bbb C}$. Our goal is to check that
 \begin{equation}
A=\prod_{k=0}^{\infty} (2a_{k}) .
\label{eq:AA}\end{equation}
Note that   under assumption \e{eq:Tr} the infinite product here  converges and $A>0$.

Let us compare the asymptotics of $\Delta (z)$ and $\Omega  (z)$ as $|z|\to\infty$. The first of them is given by \e{eq:PDhe}. So we only have to consider the Jost function $\Omega  (z)$.

    \begin{lemma}\label{HE}
    Suppose     that 
   $a_{n}=1/2$   and  $b_{n}= 0$ for all sufficiently large $n$. 
Then for all $n=-1, 0, 1,\ldots $ and $|z|\to\infty$, we have the asymptotic relation  
 \begin{equation}
f_{n}  (z)= \Big(\prod_{k=n}^\infty (2a_{k})\Big)^{-1} \z(z)^n (1+O (\z(z))).
\label{eq:HE}\end{equation}
In particular, 
 \begin{equation}
     \Omega  (z)=    A^{-1}   +O (\z(z)) .
\label{eq:HE2}\end{equation}
 \end{lemma}

 \begin{pf}
 Since $f_{n} (z)= \z^n$ for $n $ large enough, relations \e{eq:HE} are certainly satisfied for such $n$.
 Suppose that, for some $n_{0}$, relations \e{eq:HE} are true for all $n\geq n_{0}$. According to equation \e{eq:Jy} we have
  \[
 a_{n_{0}-1} f_{n_{0}-1}= \big(2^{-1} (\z+\z^{-1}) -b_{n_{0}}\big)f_{n_{0}}  -a_{n_{0} } f_{n_{0}+1}.
\]
By our assumption on $f_{n_{0}}$ and $f_{n_{0}+1}$, this implies
 \[
 a_{n_{0}-1} f_{n_{0}-1}=2^{-1}  \z^{-1}\Big(\prod_{k=n_{0}}^\infty (2a_{k})\Big)^{-1} \z^{n_{0}} (1+O (\z))
 + O (\z^{n_{0}}) 
\]
which yields \e{eq:HE} for $f_{n_{0}-1}$. In view of \e{eq:WR}, \e{eq:RV}, relation \e{eq:HE2} is a consequence of \e{eq:HE} for $n=-1$.
  \end{pf}

Thus relations \e{eq:HH} and \e{eq:AA} are verified for perturbations of finite support.
 It remains to extend them to arbitrary perturbations satisfying assumption \e{eq:Tr}.
Let us use the arguments and the notation   $H^{(N)}=H_{0}+V^{(N)}$, $\Omega^{(N)}(z) $ already used by the proof of Theorem~\ref{RRX}.  Let $\Delta^{(N)}(z)$ be the perturbation determinant for the pair $H_{0}$, $H^{(N)}$. For all $N$, we have the relation
 $$
\Delta^{(N)}(z)=   A^{(N)} \Omega^{(N)}(z)\q {\rm where}\q
 A^{(N)}=\prod_{k=0}^{N-1} (2a_{k}). 
$$
Let us pass here to the limit $N\to\infty$. Since $V^{(N)}\to V$ in the trace norm, $\Delta^{(N)}(z)\to \Delta (z)$ as $N\to\infty$. We also see  that    $\Omega^{(N)}(z)\to \Omega (z)$
by Theorem~\ref{JostP} and $A^{(N)} \to A$ under assumption \e{eq:Tr}.
This yields the desired result.

 \begin{theorem}\label{RRY}
  Under assumption \e{eq:Tr} the equality \e{eq:HH}
 is  true with constant \e{eq:AA}.
 \end{theorem}
 
 Since   $\Delta(z)$ satisfies \e{eq:PDhe}, this gives us the asymptotics of the Jost function at infinity.
 
  \begin{corollary}\label{RRB}
    Under assumption \e{eq:Tr} we have
  \[
  \lim_{| z|\to\infty}  \Omega (z)= A^{-1}.
\]
 \end{corollary}
 
 Thus, we can fix the branch of the function $\arg\Omega (z) $ for $\Im z\neq 0$ by the condition
  \[
\lim_{|z|\to\infty}\arg\Omega (z)= 0.
\]
Then in view of normalization \e{eq:PHE}, we have
 \[
\arg \Delta (z) =\arg \Omega (z),\q \Im z\neq 0.
\]

  \bigskip
  
  {\bf 5.4.}
Next, we discuss  the    spectral shift function $\xi (\lambda)$ for a pair of Jacobi operators $H_{0}$, $H$ satisfying assumption \e{eq:Tr}.
 According to Theorem~\ref{Jost} and formula 
\e{eq:HH}, the perturbation determinant $\Delta(z)$ is a continuous function of $z\in\clos\Pi$ except, possibly, the points $z=\pm 1$.  Moreover, according to \e{eq:HH5}, we have $\Omega (\lambda\pm i0)  \neq 0$.  Therefore
$\xi (\lambda)$ is also   a continuous function of $\lambda\in (-1,1)$.

A link between     $\xi (\lambda)$ and the scattering  phase $\eta(\theta)$ follows from definitions    \e{eq:AP} and \e{eq:SSF}.

\begin{theorem}\label{SzShift}
  Under assumption \e{eq:Tr},  
 the relation
  \begin{equation}
  \xi (\lambda) =\pi^{-1}\eta (\arccos\lambda) 
\label{eq:SSF3}\end{equation}
holds   for all $\lambda\in ( -1 ,\lambda)$.
 \end{theorem}

Substituting \e{eq:SSF3} into \e{eq:Sz} and taking into account \e{eq:Sz1}, we can reformulate Theorem~\ref{Sz} in terms 
of the weight function $w(\lambda)$ and   the  spectral shift function $\xi(\lambda)$.  This yields asymptotic relation \e{eq:zF}; note that    the remainder $o(1)$ in the right-hand side can be replaced by $O(\rho_{n})$.

We emphasize that $\eta(\theta)$  is a continuous function of $\theta\in (0,\pi)$, but Theorem~\ref{Jost} yields no information about its behavior as $\theta\to 0$ and $\theta\to \pi$. Comparing relations \e{eq:DS} and \e{eq:SSF3}, we however see that
 \[
\int_{0}^\pi |\eta (\theta)| \sin\theta d\theta \leq\pi  \| V\|_{1}.
\]

 \bigskip
 
      {\bf 5.5.}
   To find 
   expressions of $\tr  (  H^{n}-  H_{0}^n) $ in terms of the   spectral shift function, we only have to compare asymptotic expansions as $|z|\to\infty$ of both sides of the representation   \e{eq:DD}.   
It follows from relation
   \e{eq:PD} that
     \[
\ln \Delta (z)=-\sum_{n=1}^{\infty} n^{-1} \tr  (  H^{n}-  H_{0}^n)\, z^{-n}.
\]
Since
\[
 \int_{-\infty}^\infty \xi(\lambda) (\lambda-z)^{-1} d\lambda= -\sum_{n=1}^{\infty}  \int_{-\infty}^\infty \xi(\lambda) \lambda^{n-1} d\lambda \,  z^{-n}, 
\]
equating the coefficients at $ z^{-n}$, we see that
     \begin{equation}
  \tr  (  H^{n}-  H_{0}^n)=n   \int_{-\infty}^\infty \xi(\lambda) \lambda^{n-1} d\lambda .
\label{eq:BF1}\end{equation}

On the discrete spectrum,  the   spectral shift function  can be explicitly calculated. Indeed,
let $\lambda_{1}^{(+)}>\lambda_{2}^{(+)} >\cdots> 1$ (and $\lambda_{1}^{(-)}<\lambda_{2}^{(-)} <\cdots < - 1$) be eigenvalues of $H$ lying above the point $1$ (respectively, below the point $-1$). It follows from formula  \e{eq:ES} that
$\xi(\lambda)=n$ for $\lambda\in (\lambda_{n+1}^{(+)},\lambda_{n}^{(+)} )$  and $\xi(\lambda)=-n$ for $\lambda\in (\lambda_{n }^{(-)},\lambda_{n+1}^{(-)} )$.

Using formula  \e{eq:ES}, we find that
\[
 n\int_{1}^\infty \xi(\lambda) \lambda^{n-1} d\lambda =\sum_{k=1}^{\infty} k \big((\lambda_{k}^{(+)})^{n }-(\lambda_{k+1}^{(+)})^{n }\big)
 \]
 and, similarly, for the integral over $(-\infty,-1)$. The series here is convergent by virtue of  estimate  \e{eq:DS}. Putting together this relation with \e{eq:BF1}, we obtain the following result.
 
  \begin{theorem}\label{BF}
Let assumption   \e{eq:Tr} be satisfied. Then
\begin{equation}
  \tr  (  H^{n}-  H_{0}^n)=n   \int_{-1}^1 \xi(\lambda) \lambda^{n-1} d\lambda +\sum_{k=1}^{\infty}k \big((\lambda_{k}^{(+)})^{n }-(\lambda_{k+1}^{(+)})^{n }\big) +\sum_{k=1}^{\infty}k \big((\lambda_{k}^{(-)})^{n }-(\lambda_{k+1}^{(-)})^{n }\big).
\label{eq:BF2}\end{equation}
\end{theorem}

In view of \e{eq:SSF3}, the integral on the right can be expressed in terms of the phase function:
\[
\int_{-1}^1 \xi(\lambda) \lambda^{n-1} d\lambda=\frac{1}{\pi} \int_0^\pi \eta(\theta) \cos^{n-1}\theta \sin \theta d\lambda.
\]

      In a very general framework, formulas of type \e{eq:BF2} were studied in the book \cite{Teschl}, Chapter~6.

    \bigskip
 
   {\bf 5.6.}
   The trace formula of zero order (the Levinson theorem) requires a special discussion. Now we assume a stronger condition \e{eq:Trfd} on the coefficients of the operator $H$. Then according to Theorem~\ref{JOSF}, the corresponding perturbation determinant  $\Delta (z)$ is continuous as $z\to \pm 1$. One has to distinguish the cases  $\Delta (1)=0$ and $\Delta (-1)=0$ when the operator $H $ has threshold resonances at the points $\lambda=1$ or $\lambda=-1$. Note also (see Remark~\ref{find}) that under assumption \e{eq:Trfd} the operator $H$ has only a finite number $N$ of discrete eigenvalues.

    \begin{theorem}\label{Levinson}
Let assumption   \e{eq:Trfd} be satisfied. Then the limits $\xi(1-0)$ and $\xi(-1+0)$ exist and
\begin{equation}
 \xi(1-0)-\xi(-1+0)=N+\varkappa
\label{eq:Lev}\end{equation}
where $\varkappa=0$ if $\Delta (\pm 1)\neq 0$ for both signs, $\varkappa=1/2$ if $\Delta (\pm 1)=0$ for one of the signs and $\varkappa=1$ if $\Delta (\pm 1)=0$ for both   signs.
\end{theorem}

    \begin{pf} 
    Let us consider a contour which consists of a part 
\[
{\cal C}_{r_{- }, r_{+}}=[-1+r_{- }+i0, 1-r_{+}+i0]\cup [-1+r_{-}-i0, 1-r_{+}-i0]\cup \{|z-1]=r_{+}\}\cup \{|z+1]=r_{-}\}
\]
 encircling the cut along $[-1,1]$ and of the circle $|z|=R$. We suppose that $r_{-}$ and $r_{+}$ are sufficiently small,  $R$ is sufficiently large, and we go   in the clockwise over  ${\cal C}_{r_{- }, r_{+}}$ and   in the counter-clockwise direction over  $|z|=R$. Since $\Delta(z)$ has simple zeros at the eigenvalues of $H$,  the argument principle implies that
 \begin{equation}
 \var_{{\cal C}_{r_{- }, r_{+}}} \arg\Delta (z)+ \var_{|z|=R} \arg\Delta (z)=2\pi N .
\label{eq:Lev1}\end{equation}
According to \e{eq:PHE}  we have $\var_{|z|=R} \arg\Delta (z)\to 0$ as $R\to\infty$. It follows from  \e{eq:SSF} that
 \begin{equation}
 \var_{{\cal C}_{r_{- }, r_{+}}} \arg\Delta (z) =2\pi (\xi (1-r_{+})- \xi (1+r_{-}))+ \var_{ |z-1|=r_{+}} \arg\Delta (z) + \var_{ |z+1|=r_{-}}  \arg\Delta (z).
\label{eq:Lev2}\end{equation}

Let us now pass here to the limit $r_{-}, r_{+}\to 0$. If $\Delta (\pm 1)\neq 0$, then $ \var_{ |z\mp 1|=r_{\pm}} \arg\Delta (z)\to 0$.   In the case $\Delta (\pm 1)=0$ it follows from Theorem~\ref{resin1} that
\[
\lim_{r_{\pm} \to 0} \var_{ |z\mp 1|=r_{\pm}} \arg\Delta (z) =\pi.
\]
Thus relation \e{eq:Lev} is a direct consequence of \e{eq:Lev1} and \e{eq:Lev2}. These arguments prove also the existence of the limits $\xi(1-0)$ and $\xi(-1+0)$.
    \end{pf}

 \section{Scattering theory}

 {\bf 6.1.} 
First,  we briefly recall basic notions of scattering theory.
 We refer, for example,  to the book \cite{Ya} for a more    complete presentation of this  material.
 
 The wave operators $W_{\pm} (H , H_{0})$   for a pair of self-adjoint operators $H_{0}$, $H $ in a Hilbert space ${\cal H}$  are defined as strong limits
   \begin{equation}
  W_{\pm} (H , H_{0})=\slim_{t\to \pm\infty}e^{iH t}e^{-iH_{0}t};
\label{eq:WO}\end{equation}
here $H_{0}$ is supposed to be absolutely continuous.
 Under the assumption of the existence of limits \e{eq:WO}, $ W_{\pm} (H , H_0)$ are isometric operators and enjoy the intertwining property
 $H W_{\pm} (H , H_{0})=W_{\pm} (H , H_0) H_{0}$.   The wave operator $W_{\pm} (H , H_{0})$  is called complete  if its range coincides with the absolutely  continuous subspace of the operator $H $. The  scattering operator
 \begin{equation}
 {\bf S} =W_{+} (H , H_{0})^* W_{-} (H , H_{0})
 \label{eq:WOS}\end{equation}
 commutes with $H_{0}$,  ${\bf S} H_0 =H_0{\bf S} $,  and it is unitary if both  wave operators $W_{+} (H , H_{0})$  
 and $W_- (H , H_{0})$  are complete. 
 
 To define the scattering matrix, we suppose for definiteness  that the spectrum of $H_{0}$ is simple and coincides with $[-1,1]$. Let $F_{0}: {\cal H} \to L^2 (-1,1)$ be unitary and $F_{0} H= {\sf A}F_{0}$ where ${\sf A}$ is the operator of multiplication by $\lambda$ in the space $L^2 (-1,1)$. Then  the scattering matrix
 $S (\lambda )\in{\Bbb C}$ is defined by the relation
   \begin{equation}
  (F_{0} {\bf S}  f)(\lambda)=S  (\lambda ) (F_0 f)(\lambda),\q \lambda\in (-1,1);
   \label{eq:SM}\end{equation}
obviously,   $|S  (\lambda )|=1$ if $ {\bf S}$ is unitary.  
    Note that the scattering matrix does not depend on the diagonalization $F_{0}$ of $H_0$ if it has simple spectrum.

Let us come back to Jacobi operators. Under condition \e{eq:Tr} which is always assumed in this section,
 the perturbation $V=H-H_{0}$  is trace class.  Therefore the wave operators $W_{\pm} (H , H_{0})$ exist and are  complete by the classical Kato-Rosenblum theorem.
 
 Our goal here is to obtain representations of the wave operators and of the scattering matrix in terms of the polynomials $P_{n}(\lambda)$ and of the Jost function $\Omega(\lambda+i0)$ for $\lambda\in (-1,1)$.
  Such expressions are known as stationary representations.   As a by-product of our considerations, we will also give a direct proof  of the existence and completeness of the wave operators.
 
\bigskip

 {\bf 6.2.}
 Let $d\rho(\lambda)$ be the spectral measure of the operator $H$. We  define a mapping $U: \ell^2 ({\Bbb Z}_{+})\to L^2 ({\Bbb R}; d\rho)$ by the formula 
 \[
 (Ue_{n})(\lambda)=P_{n} (\lambda).
 \]
This mapping is isometric according to \e{eq:J5}. It is also unitary because  the set of all polynomials  $P_n(\lambda) $, $n\in {\Bbb Z}_{+}$,  is dense in $L^2 ({\Bbb R}; d\rho)$. Finally, the intertwining property 
 \begin{equation}
(U H f) (\lambda)= \lambda (U f)(\lambda)
\label{eq:IP}\end{equation}
  holds. Indeed, it suffices to check it for $f=e_{n}$ when according to definition \e{eq:J} of the operator $H$, $(U H e_{n}) (\lambda)$ coincides with the left-hand side of \e{eq:Py}  (where $z=\lambda$) while $\lambda (U e_{n})(\lambda)$ equals its right-hand side.

Next, we reduce the absolutely continuous part of the operator $H$ to the operator  ${\sf A}$ of multiplication by $\lambda$ in $L^2 (-1,1 )$.   To that end, we put
   \begin{equation}
\psi_{n}(\lambda)= \sqrt{w(\lambda)}  P_{n} (\lambda), \q \lambda\in (-1,1),
\label{eq:UF}\end{equation}
where $w(\lambda)$ is the weight function     defined by \e{eq:SF1} and introduce a mapping $F : \ell ^2 ({\Bbb Z}_{+})\to L^2 (-1,1 )$  by the formula
 \begin{equation}
 (F  e_{n})(\lambda)=\psi_{n} (\lambda), \q \lambda\in (-1,1).
\label{eq:psX}\end{equation}
The operator $F^* :  L^2 (-1,1)\to  \ell^2 ({\Bbb Z}_{+})$ adjoint to $F $  is given by the equality
   \begin{equation}
(F ^* g)_{n}=   \int_{-1}^1  \psi_{n} (\lambda ) g(\lambda) d\lambda,\q n\in {\Bbb Z}_{+}.
\label{eq:UFf}\end{equation}
 According to \e{eq:UF} it follows from the unitarity of the operator $U $  that
  \begin{equation}
F  F ^*=I,\q  F ^* F   =  E (-1,1) ,
\label{eq:ps3}\end{equation}
where $E (-1,1)$ is the spectral projection of the operator $H $ corresponding to the interval $(-1,1)$.
Thus the operator $F^*$ is isometric, and it is   unitary if $H$ has no point spectrum.  
By virtue of the  equation \e{eq:IP} the intertwining property
  \begin{equation}
H F ^*=F ^* {\sf A}
\label{eq:inw1}\end{equation}
holds. 

For the free operator  $H_0$,  the corresponding unitary mapping  $F_{0} : \ell ^2 ({\Bbb Z}_{+})\to L^2 (-1,1 )$ is defined by the formula
  \begin{equation}
 (F_{0}  e_{n})(\lambda)=\psi_{n}^{(0)} (\lambda)\q {\rm where} \q \psi_{n}^{(0)} (\lambda)=\sqrt{w_{0}(\lambda)} P_{n}^{(0)} (\lambda), \q \lambda\in (-1,1).
\label{eq:pX}\end{equation}
Here $w_{0} (\lambda) $ is defined by \e{eq:ww}
and $P_{n}^{(0)}  (\lambda)$ are normalized Chebyshev polynomials of the second kind.
 
\bigskip

  {\bf 6.3.} 
  The presentation below is very close to Section 4.2 of the book \cite{YA}  where the continuous case was considered.
We start with     two   elementary facts. Recall that  $\z (\lambda\pm i0)$ and $\rho_{n}$ are defined by formulas  \e{eq:zt} and  \e{eq:sig}, respectively.
 
  \begin{lemma}\label{CheW}
 Let $\int_{-1}^1 | g(\lambda) |^2 \sqrt{1-\lambda^2}d \lambda<\infty$. Then
     \begin{equation}
\lim_{t\to\pm\infty}\sum_{n=0}^\infty \big| \int_{-1}^1 \z (\lambda\mp i0)^n e^{-i\lambda t} g(\lambda) d\lambda \big|^2
=0.
\label{eq:opm1}\end{equation}
  \end{lemma}

  \begin{pf}
  It is more convenient to work in the variable $\theta$ using that
    \begin{equation}
\int_{-1}^1\z (\lambda\mp i0)^n e^{-i\lambda t} g(\lambda) d\lambda= \int_0^\pi e^{\pm i n\theta}e^{-i t\cos\theta } g(\cos\theta) \sin\theta  d\theta.
\label{eq:theta}\end{equation}
  The Parseval identity implies the uniform in $t$ estimate
       \begin{multline*}
\sum_{n=0}^\infty \big| \int_{-1}^1\z (\lambda\mp i0)^n e^{-i\lambda t} g(\lambda) d\lambda \big|^2
\\
\leq 2\pi
\int_0^\pi | g(\cos\theta) |^2\sin^2\theta d \theta = 2\pi
\int_{-1}^1 | g(\lambda) |^2 \sqrt{1-\lambda^2}d \lambda  .
\end{multline*}
Therefore it suffices to check \e{eq:opm1} for $g\in C_{0}^\infty ((-1,1))$.

Integrating in \e{eq:theta} by parts and observing that
\[
\big| \pm n+ t \sin\theta \big|\geq c \big|n+ |t| \big|  \q {\rm for}  \q \pm t>0
\]
where $\theta\in [\varepsilon_{0}, \pi-\varepsilon_{0}]$, $\varepsilon_{0}>0$ and   $c=c (\varepsilon_{0})>0$, we find an estimate
  \[
\big| \int_0^\pi e^{\pm i n\theta}e^{-i t\cos\theta } g(\cos\theta) \sin\theta  d\theta \big|
\leq C (n+ |t|)^{-1},\q n\in {\Bbb Z}_{+},  \q \pm t>0.
\]
    It follows  that
 \[
 \sum_{n=0}^\infty \big| \int_{-1}^1 \z (\lambda\mp i0)^n e^{-i\lambda t} g(\lambda) d\lambda \big|^2
\leq C  \sum_{n=0}^\infty (n+ |t|)^{-2}\leq C_{1} |t|^{-1},  \q \pm t>0,
\]
which implies \e{eq:opm1}.
      \end{pf}

       \begin{lemma}\label{CheWz}
         Under assumption \e{eq:Tr},  suppose additionly that
            \begin{equation}
\sum_{n=0}^\infty \rho_{n}^2 <\infty.
\label{eq:rho}\end{equation}
  Let  $g\in L^1 (-1,1)$ and $\supp g \subset (-1,1)$. Then
     \begin{equation}
\lim_{|t |\to \infty}\sum_{n=0}^\infty \big| \int_{-1}^1\big(f_{n}(\lambda\pm i0)- \z (\lambda\pm i0)^n \big) e^{-i\lambda t} g(\lambda) d\lambda \big|^2
=0.
\label{eq:opm2}\end{equation}
  \end{lemma}
  
   \begin{pf}
    Theorem~\ref{Jost}  yields the estimate
   \[
   \big| \int_{-1}^1\big(f_{n}(\lambda\pm i0)- \z (\lambda\pm i0)^n \big) e^{-i\lambda t} g(\lambda) d\lambda \big|\leq C \rho_{n}
   \]
   where $C$ does not depend on $t$. For every $n$, the integrals in \e{eq:opm2} tend to zero as $|t|\to\infty$ by the Riemann-Lebesgue lemma. So, under assumption \e{eq:rho}, the dominated convergence theorem implies \e{eq:opm2}.
      \end{pf}

    Let us now set  
     \begin{equation}
 \sigma_\pm(\lambda) = \frac{\Omega(\lambda\pm i0)} { |\Omega(\lambda\pm i0)|}=e^{\pm i\eta(\arccos \lambda)},\q \lambda\in (-1,1) .  
 \label{eq:CH1r}\end{equation} 
The operator $ \Sigma_\pm(\lambda)$ of multiplication by $ \sigma_\pm(\lambda)$ is of course unitary in the space $ L^2 (-1,1)$.  
 
 Using definitions \e{eq:SF1} and   \e{eq:UF} of the functions $w  (\lambda)$  and 
  $\psi_{n}  (\lambda)$,  we rewrite relation   \e{eq:HH4} as
      \begin{equation}
\psi_{n}  (\lambda)=\frac{1} {i \sqrt{2\pi} \sqrt[4]{1-\lambda^2} }\Big(-\z (\lambda+i0) \sigma_{-} (\lambda)f_{n} (\lambda+i0) + \z (\lambda- i0)\sigma_+(\lambda)f_{n} (\lambda-i0) \Big).  
 \label{eq:YY}\end{equation}
  For the operator $H_{0}$ and the function $\psi_{n}^{(0)}  (\lambda)$   defined in \e{eq:pX}, this reduces to the formula
       \begin{equation}
\psi_{n}^{(0)}  (\lambda)=\frac{1} {i \sqrt{2\pi} \sqrt[4]{1-\lambda^2} }\Big( -  \z(\lambda+i0)^{n+1} + \z(\lambda- i0)^{n+1}  \Big)  .
 \label{eq:YY1}\end{equation}

  Recall that the operators $F$ and  $F_{0}$  are defined  by relations  \e{eq:psX} and \e{eq:pX}.  
  Then  $F^*$ is  given by relation  \e{eq:UFf} and similarly for  $F^*_{0}$. The operator ${\sf A}$ acts as multiplication by $\lambda$ in the space $  L^2 (-1,1)$.

      \begin{lemma}\label{CheW1}
For all $g\in L^2 (-1,1)$, we have
    \begin{equation}
 \lim_{t \to\pm\infty} \| (F^* \Sigma_{\pm}   -F_0^*)e^{-i{\sf A}t}g\|=0.
 \label{eq:CH1rr}\end{equation}
  \end{lemma}
  
    \begin{pf}
    Since $\sigma_{+} (\lambda) \sigma_- (\lambda) =1$,  it follows from relations  \e{eq:YY}  and  \e{eq:YY1}  that, for both signs  $``\pm"$,
      \begin{multline}
i \sqrt{2\pi} \sqrt[4]{1-\lambda^2} \big( \psi_{n}  (\lambda) \sigma_\pm(\lambda)-\psi_{n}^{(0)}  (\lambda)\big)\\
=    
 \pm ( \sigma_{\pm} (\lambda)^2 -1)  \z(\lambda\mp i0)^{n+1} 
 +\sigma_{\pm} (\lambda)     r_{n}    (\lambda) 
 \label{eq:YY2}\end{multline}
 where $ r_{n}    (\lambda)= -    r_{n}^{(+)}   (\lambda)
+    r_{n}^{(-)}   (\lambda)$ and
      \[
r_{n}^{(\pm)}   (\lambda)= \z(\lambda\pm i0)\sigma_{\pm} (\lambda) \big( f_{n} (\lambda\pm i0) -\z(\lambda\pm i0)^{n }\big)    .
   \]
According to \e{eq:YY2} we have
   \begin{multline}
(F^* \Sigma_\pm  -F_0^*)e^{-i{\sf A}t}g)_{n } =  \pm \int_{-1}^1 ( \sigma_{\pm} (\lambda)^2 -1) \z(\lambda\mp i0)^{n+1} e^{-i\lambda t}
\tilde{g} (\lambda)d\lambda
\\
+ \int_{-1}^1 r_{n} (\lambda)   \sigma_\pm (\lambda) e^{-i\lambda t}
\tilde{g} (\lambda)d\lambda
 \label{eq:YX}\end{multline}
where  
\[
\tilde{g} (\lambda)=\frac{g(\lambda)}{i \sqrt{2\pi} \sqrt[4]{1-\lambda^2}} .
\]

Let us come back to relation  \e{eq:CH1rr} where we may assume that $\supp g \subset (-1,1)$.
If $t\to+\infty$ (if $t\to -\infty$), we use \e{eq:YX} for the upper (lower) sign.   Then the contribution to the norm in $\ell^2 ({\Bbb Z}_{+})$ of the first   term in the right-hand side of 
\e{eq:YX} tends to zero according to Lemma~\ref{CheW}. The contribution  of the second   term tends to zero according to Lemma~\ref{CheWz}.
  \end{pf}

   \begin{theorem}\label{CheW2}
     Let assumptions \e{eq:Tr} and  \e{eq:rho} be satisfied.  
Then the strong limits \e{eq:WO}  exist and
     \begin{equation}
W_{\pm}(H ,H_{0})= F ^* \Sigma_{\pm}    F_{0}.
 \label{eq:CH1pr}\end{equation}
  \end{theorem}
  
     \begin{pf}
   We have to check that
          \[
 \lim_{t \to\pm\infty} \| e^{iH t} e^{-iH_0 t}f - F ^* \Sigma_{\pm}   F_0 f\|=0
\]
 for all $f\in\ell^2({\Bbb Z}_{+})$. In view of the intertwining property \e{eq:inw1} this relation can be rewritten as
   \begin{equation}
 \lim_{t \to\pm\infty} \| e^{-iH_0 t}f - F^* \Sigma_{\pm}    e^{-i{\sf A}t} F_{0}f\|=0.
 \label{eq:wo1}\end{equation}
 Set now $g=F_{0}f$. Then again in view of the intertwining property \e{eq:inw1} for   $H_{0}$, we see that relations
 \e{eq:CH1rr} and \e{eq:wo1}  are equivalent.
          \end{pf}

  It follows from relations \e{eq:ps3} and  \e{eq:CH1pr}  that the scattering operator \e{eq:WOS} is given by the equality
  \[
  {\bf S} = F_0^* \Sigma_{+} ^{*} \Sigma_{-}   F_{0}.
  \]
  Putting together this relation with the definition  \e{eq:SM} of  the scattering matrix and definition  \e{eq:CH1r} of $\sigma_\pm (\lambda )$, we can state the following result. Recall also that the scattering phase $\eta (\theta)$ was defined by formula \e{eq:AP}.
  
   \begin{theorem}\label{WS}
     Under assumptions of Theorem~\ref{CheW2},   the scattering matrix for the pair $H_{0}$, $H$ satisfies the equality
    \begin{equation}
S  (\lambda)=  \frac{\Omega (\lambda- i0)} {\Omega (\lambda+ i 0)}=e^{-2 i \eta (\arccos\lambda)},\q \lambda\in (-1,1) .
 \label{eq:SS}\end{equation}
  \end{theorem}
  
In view of relation  \e{eq:SSF3}, formula \e{eq:SS} can be rewritten as the Birman-Kre\u{\i}n formula
\[
S  (\lambda)=   e^{-2 i \pi \xi( \lambda)}.
\]
  
  \bigskip
  
   {\bf 6.4.}
    Theorems~\ref{CheW2} and \ref{WS} remain true without additional assumption \e{eq:rho}. The proof of this fact requires some tools of abstract scattering theory. We give only basic ideas of one of its possible proofs.

    The following result is a direct consequence of Theorem~\ref{ORP}.
    
    \begin{lemma}\label{H-S}
Let    an operator $Q$ in the space $\ell^2 ({\Bbb Z}_{+})$ be defined by the equality  $(Q  f)_{n}= q_{n} f_{n}$ where $q_{n}=  \bar{q}_{n}$ and $\sum_{n=0}^\infty q_{n}^2 <\infty$. Under assumption \e{eq:Tr} the operator-valued function  $Q  R  (z) Q $ defined for $\Im z\neq 0$ is continuous in the Hilbert-Schmidt norm up to the cut along $[-1,1]$ with possible exception of the  points $z=\pm 1$. In particular, the operator $Q$ is   $H$-smooth in the sense of Kato on every compact subinterval of $(-1,1)$.
  \end{lemma}
  
  Under assumption \e{eq:Tr} Lemma~\ref{H-S} is true with
  \[
  q_{n}= (| a_{n}-1/2|+|b_{n}|)^{1/2}.
  \]
  Note that the perturbation $V=Q G  Q $  where $G $ is a bounded  operator in the space $\ell^2 ({\Bbb Z}_{+})$.
   Using Theorems~5.3.4 and 5.6.1  in \cite{Ya}, we can now deduce from  Lemma~\ref{H-S} the following assertion.

 \begin{theorem}\label{CA}
 The wave operators $W_{\pm}  (H,H_{0})$ exist, are complete and representations \e{eq:CH1pr}, \e{eq:SS} are true under the only assumption
    \e{eq:Tr}.
  \end{theorem}

    \section{The Szeg\"o function and the Case sum rules}
 
 Here we establish a link between    the perturbation determinant $\Delta(z)$ and the Szeg\"o function $D(\z)$. 
 We also show that the Case sum rules are direct consequences of this link.
 Since the Jost function $\Omega(z)$ is related to $\Delta(z)$ by simple formula \e{eq:HH}, we do not discuss it here.
 As usual, $|\z|<1$ and $2z=\z+\z^{-1}$.
 
 \medskip
 
    {\bf 7.1.}
    We define 
  the Szeg\"o function $D(\z)$   by the formula
       \begin{equation}
D(\z)=  \exp\Big(\frac{1 }{4\pi}\int_{-\pi}^\pi  \frac{  e^{i\theta}+\z }{  e^{i\theta}-\z } \ln \big( w(\cos\theta)|\sin\theta|\big)d\theta\Big),\q |\z|< 1,
\label{eq:SzY}\end{equation}
where $w(\lambda)$ is the weight function  \e{eq:SFx}. This is exactly formula (10.2.10) (see also Theorem~12.1.2) in \cite{Sz}, but in contrast to \cite{Sz} we do not suppose that $H$ has no eigenvalues. Of course definition \e{eq:SzY} requires that
\[
\int_{-\pi}^\pi | \ln \big( w(\cos\theta)|\sin\theta|\big)|d\theta<\infty
\]
or, equivalently, \e{eq:Szeg} be satisfied. This is known as the Szeg\"o condition.

     Recall the standard Jensen-Poisson representation of functions $f(\z)$ analytic in the unit disc $|\z| <1$. If $f(\z)$ is continuous in the closed disc $|\z|\leq 1$ and $\Im f(0)=0$, then
       \begin{equation}
f(\z)=  \frac{1 }{2\pi}\int_{-\pi}^\pi  \frac{  e^{i\theta}+\z }{  e^{i\theta}-\z }   \Re f(e^{i\theta})d\theta .
   \label{eq:JP}\end{equation}
 In particular, we see that  
   \begin{equation}
   \begin{split}
    \ln   (1+ \z) &=
 \frac{1 }{2\pi}\int_{-\pi}^\pi  \frac{  e^{i\theta}+\z }{  e^{i\theta}-\z } \ln \cos (\theta/2) d\theta, 
\\
    \ln   (1- \z) &=
 \frac{1 }{2\pi}\int_{-\pi}^\pi  \frac{  e^{i\theta}+\z }{  e^{i\theta}-\z } \ln| \sin (\theta/2) |d\theta. 
    \end{split}
\label{eq:wwz}\end{equation}
We here fix $\arg(1+\z)$ and $\arg(1- \z)$  by the condition $\arg 1=0$.

   Set $\wt\Delta(\z)=\Delta(z)$. It follows from Corollary~\ref{JOST} that the function $\wt\Delta(\z)$ is analytic in the unit disc and is continuous up to the unit circle  with a possible exception of the points $\pm 1$.
Moreover, according to  \e{eq:HH5},  $\wt\Delta(\z)\neq 0$ if $|\z| =1$ but $\z\neq \pm 1$.  Under assumption \e{eq:Tr}, an information on the behavior of $\wt\Delta(\z)$ as  $\z\to\pm 1$ follows from Theorem~9.14 in 
 \cite{KS} where  it is shown  that the function $\wt\Delta(\z)$ belongs to the Nevanlinna class $N$ and does not have a singular inner component. We refer to the book \cite{Duren}, Chapter~2, for precise definitions of these notions.
  We also note that $\wt\Delta(0)=1$ according to \e{eq:PDhe}.

      Let $\lambda_{k}$ be eigenvalues (lying on $(-\infty, -1)\cup(1,\infty)$) of the operator $H$. In contrast to Subsection~5.5, we do not distinguish positive and negative eigenvalues in notation and suppose that
      $|\lambda_{1}|\geq |\lambda_{2}|  \geq \cdots >1$. 
          The  numbers $\mu_{k}\in (-1,1)$ defined by $\mu_{k}+ \mu_{k}^{-1} =2\lambda_{k}$ are zeros of the function $\wt\Delta(\z)$. 
      Since $\wt\Delta\in N $, we have
     \begin{equation}
\sum_{k=1}^\infty (1-|\mu_{k}|)<\infty,
\label{eq:Bla1}\end{equation}
so that  the   Blaschke product
      \begin{equation}
B(\z) =\prod_{k=1}^\infty\frac{\mu_k}{|\mu_{k}|}\frac{\mu_{k}-\z}{1-\mu_k\z},\q |\z|<1,
\label{eq:Bla}\end{equation}
 is well defined. The function $B(\z)$ is continuous in the closed disc $|\z|\leq 1$ except, possibly, the points $\pm 1$ and 
 $|B(\z)|=1$ for $|\z| =1$. Note that the proof   of the inclusion $\wt{\Delta}\in N$ in \cite{KS} relied on the condition \e{eq:Bla1}   established under assumption \e{eq:Tr} earlier in \cite{Hu-S}.

 Let us define the outer function
      \begin{equation}
G(\z)= \exp\Big(\frac{1 }{2\pi}\int_{-\pi}^\pi  \frac{  e^{i\theta}+\z }{  e^{i\theta}-\z } \ln | \wt{\Delta}(e^{i\theta})|d\theta\Big)
   \label{eq:Bla2}\end{equation}
   where the function $\ln | \wt{\Delta}(e^{i\theta})|$ belongs to $L^1 (-\pi,\pi)$ 
   because $\wt{\Delta}\in N$. This is equivalent to the Szeg\"o condition \e{eq:Szeg}
   since, by relations \e{eq:SF1}  and \e{eq:HH},  
\[
| \wt{\Delta}(e^{i\theta})|^{2}= | {\Delta}(\cos\theta)|^{2}= A^{2}   \frac{2}{\pi}\frac{|\sin\theta|}{w(\cos\theta)},
\]
whence
 \begin{equation}
2 \ln | \wt{\Delta}(e^{i\theta})|=-  \ln (   w(\cos\theta)|\sin\theta|)  +\ln (2\pi^{-1} A^2 \sin^2\theta).
 \label{eq:wj}\end{equation}

 It follows from Theorem~2.9 in \cite{Duren} that the factorization
    \begin{equation}
  \wt{\Delta}(\z)=      
B(\z)    G(\z)  
 \label{eq:fact}\end{equation}
holds. Of course \e{eq:Bla2} is the Jensen-Poisson representation \e{eq:JP} of the function $ f(\z)= \ln \big(\wt{\Delta}(\z)/      
B(\z) \big)$, but, since this function is not continuous in the disc $|\z|\leq  1$,  its justification relies on the properties of $\wt{\Delta} $ stated above.

Let us now compare definitions    \e{eq:SzY} and \e{eq:Bla2}.  
We substitute   expression \e{eq:wj} into \e{eq:SzY} and take into account    formulas \e{eq:wwz}. Thus, using factorization    \e{eq:fact}, we arrive at the following result.

 \begin{theorem}\label{PertSze}
Let assumption   \e{eq:Tr} be satisfied, and let $|\z |< 1$. Set $\wt\Delta(\z)=\Delta(2^{-1}(\z+\z^{-1}))$ where $\Delta(z)$ is the perturbation determinant  \e{eq:Tr1}. Define the Blaschke product  $B(\z)$ by formula  \e{eq:Bla}, the Szeg\"o function $D(\z)$ --  by \e{eq:SzY}, and the product $A$  --  by \e{eq:AA}.
Then the identity
  \begin{equation}
\wt{\Delta} (\z)= A B(\z)\frac{1-\z^2}{\sqrt{ 2\pi} D(\z)}
\label{eq:SzD3}\end{equation}
is true.
\end{theorem}

We emphasize that Theorem~\ref{PertSze} is a direct consequence of classical results on the factorization of functions in the Nevanlinna class combined with the analytical results of \cite{KS}.

According to \e{eq:SzD3} formulas \e{eq:DD} and \e{eq:SzY}    provide two different representations of an essentially the same object named the perturbation determinant $\Delta(z)$   or  the Szeg\"o function $D(\z)$. In view of \e{eq:SSF} the first of them is given 
  in terms of $\arg \Delta(\lambda+i0)$ while according to \e{eq:SF1} and \e{eq:HH} the second representation is stated  
  in terms of  $\ln |\Delta(\lambda+i0)|$. Obviously, these two functions are  harmonic conjugate.

It is of course possible to  rewrite  representation  \e{eq:SzY} in terms of the variable $z=2^{-1}(\z+\z^{-1})\in\Pi$. Let us also introduce   the weight function  $w_{0}(\lambda)$ of the operator $H_{0}$; it is given by formula \e{eq:ww}.
  Taking  into account that the function $w(\cos\theta)|\sin\theta|$ is even, making the change of variables $\lambda=\cos \theta$ and using formulas  \e{eq:wwz}, we see that
    \begin{equation}
D(\z)= \frac{1-\z^2}{\sqrt{ 2\pi}  }\exp\Big(\frac{1-\z^2}{2\pi}\int_{-1}^1  \frac{\ln \big(w(\lambda)/w_{0}(\lambda)\big)}{1-2\z \lambda+ \z^2}\frac{d\lambda}{\sqrt{1-\lambda^2}}\Big)
\label{eq:SzD1}\end{equation}
or, equivalently,
   \[
D(\z (z))= \z (z) \sqrt{\frac{2(z^{2}-1)}{  \pi}  }\exp\Big(-\frac{\sqrt{z^2-1}}{2\pi}\int_{-1}^1  \frac{\ln \big(w(\lambda)/w_{0}(\lambda)\big)}{ \lambda-z}\frac{d\lambda}{\sqrt{1-\lambda^2}}\Big).
\]

\bigskip
 
   {\bf 7.2.}
   Let us, finally, obtain the Case sum rules for the pair $H_{0}$, $H$.
  Putting together \e{eq:SzD3} and \e{eq:SzD1}, we see that
   \begin{equation}
\ln  \wt\Delta (\z) -\ln B(\z) -\ln A= \frac{\z^2-1}{2\pi}\int_{-1}^1  \frac{\ln \big(w(\lambda)/w_{0}(\lambda)\big)}{1-2\z \lambda+ \z^2}\frac{d\lambda}{\sqrt{1-\lambda^2}} .
\label{eq:trace}\end{equation}
Let us   compare the behavior of both sides of this equation as $\z\to 0$.  Since
$ \wt\Delta (0) =1$,  setting $\z=0$ and using definition \e{eq:Bla}, we obtain the identity
  \begin{equation}
 \sum_{k=1}^\infty \ln (2a_{k}) + \sum_{k=1}^\infty \ln |\mu_{k}| =  \frac{1}{2\pi}\int_{-1}^1   \ln \big(w(\lambda)/w_{0}(\lambda)\big) \frac{d\lambda}{\sqrt{1-\lambda^2}} .
\label{eq:case}\end{equation}
This Case sum rule of zero order is of course different from the Levinson theorem \e{eq:Lev}.

More generally, let us consider the asymptotic expansions of both sides of \e{eq:trace} as $\z\to 0$ and compare the coefficients at the same powers of $\z$. According to Theorem~2.13 in \cite{KS} we have
\[
\ln  \wt\Delta (\z) =-2 \sum_{n=1}^\infty n^{-1} \tr \big( T_{n} (H)- T_{n} (H_{0})\big) \z^n 
\]
where $T_{n}(\lambda)=\cos (n\arccos\lambda)$ are the Chebyshev polynomials of the first kind.
It directly follows from definition \e{eq:Bla} that
     \[
\ln B(\z) =\sum_{k=1}^\infty \ln |\mu_k| + \sum_{n=1}^{\infty}n^{-1}\big(\sum_{k=1}^{\infty} (\mu_k^{n}-\mu_k^{-n})\big)\z^{n}
\]
where the series over $k $ are convergent due to the condition  \e{eq:Bla1}.
Finally, we use formula (10.11.29) in \cite{BE}: 
  \[
  \frac{1-\z^2} {1-2\z \lambda+ \z^2} =1+2 \sum_{n=1}^\infty T_{n} (\lambda) \z^n. 
\]
Thus the equality of the coefficients at $\z^{n}$ in the left-  and right-hand sides of \e{eq:trace} yields the identity
\begin{multline}
\tr \big( T_{n} (H)- T_{n} (H_{0})\big) =-\frac{1}{2} \sum_{k=1}^\infty  (\mu_k^{n}-\mu_k^{-n})
\\+ \frac{n}{2\pi}
\int_{-1}^1   \ln \big(w(\lambda)/w_{0}(\lambda)\big) T_{n} (\lambda) \frac{d\lambda}{\sqrt{1-\lambda^2}},\q n=1,2,\ldots. 
\label{eq:case1}\end{multline}

The  trace identities \e{eq:case} and \e{eq:case1} known as the Case sum rules are not new. They    were obtained by him  in  \cite{Case} and rigorously proven in \cite{KS}.
We note that, in the paper \cite{KS},  the identities \e{eq:case} and \e{eq:case1}  were first checked for finite rank perturbations  $H-H_{0}$ and then  \e{eq:case} and \e{eq:case1} (for $n=1$) were used for
  the proof of the inclusion $\wt{\Delta}\in N$.

\bigskip
 
    {\bf 7.3.} 
    The Szeg\"o   condition \e{eq:Szeg} implies that the weight function $w(\lambda)$ does not vanish too rapidly as $\lambda\to -1$ and $\lambda\to 1$ or, to put it differently,  $|\Omega( \lambda+i0)|$   does not tend to infinity too rapidly.
  The example of the Pollaczek polynomials shows that condition \e{eq:Szeg}  may be violated without assumption \e{eq:Tr}.
  
  Recall that the normalized Pollaczek polynomials are defined (see, e.g., Appendix to \cite{Sz}) by recurrent relations \e{eq:Py} with
     \begin{equation}
 a_{n}=\frac{n+1}{\sqrt{(2n+2\alpha +1)(2n+2\alpha +3)}},\q  b_{n}=-\frac{2\beta}{ 2n+2\alpha +1 }; 
\label{eq:Poll1}\end{equation}
here the parameters $\alpha, \beta\in {\Bbb R}$ and $\alpha> |\beta|$.
It follows   that
    \begin{equation}
 a_{n}=2^{-1} -\alpha (2n)^{-1}+ O (n^{-2}),\q  b_{n}=-\beta n^{-1}  +  O (n^{-2}).
\label{eq:PolX}\end{equation}
The corresponding  normalized  weight function is given by the formula
 \[
w(\lambda)= (\alpha+ 1/2) e^{(2\theta-\pi) h(\theta)}\big( \cosh (\pi h(\theta) \big)^{-1}\q {\rm where}\q  h(\theta) =(\alpha \cos\theta +\beta) (\sin\theta)^{-1}
\]
and as usual $\lambda=\cos\theta$.
It is easy to see that
  \begin{equation}
\ln w(\lambda)= - \pi (\alpha+\beta)\theta^{-1}+ O(1) 
\label{eq:PolCou}\end{equation}
as $\lambda\to1-0$ and a similar formula is true as $\lambda\to -1+ 0$.

It follows from \e{eq:PolX} that $V=H-H_{0}$ is Hilbert-Schmidt, but the series $\sum_{n} (a_{n}-1/2)$ and $\sum_{n} b_{n}$ are divergent; in particular, assumption \e{eq:Tr} is not satisfied.  
According to \e{eq:PolCou} the Szeg\"o condition \e{eq:Szeg} is violated for Pollaczek polynomials. This is consistent with the classical theorem of Szeg\"o, Shohat, Geronimus, Kre\u{\i}n  and Kolmogorov; see, e.g., Theorem~4 in \cite{KS}. On the other hand, relation \e{eq:PolCou}  implies that
 \[
\int_{-1}^1 \ln w(\lambda) (1-\lambda^2)^{1/2} d\lambda> -\infty. 
\]
This is consistent with  Theorem~1 in \cite{KS}.

\bigskip
 
   {\bf 7.4.}
   For the continuous operator ${\bf H}= -D^2 + b(x)$,  an analogue of the Szeg\"o   condition means that the corresponding weight function 
${\bf w}(\lambda) $ (for its definition, see, e.g., subsection~6 of $\S$ 4.1 in \cite{YA}) does not vanish too rapidly as $\lambda\to +0$, that is,
 \begin{equation}
\int_{0}^1 \ln {\bf w}(\lambda) \lambda^{-\alpha} d\lambda> -\infty 
\label{eq:Szeg1}\end{equation}
where $\alpha=1/2$ if $b\in L^1 ({\Bbb R}_{+})$. Note however that we were unable to find this assertion in the literature.

Condition \e{eq:Szeg1} appears to be sufficiently sharp. Indeed, suppose that $b$ is a non-negative smooth function such that 
\[
b(x)=b_{0} x^{-\rho}, \q b_{0} >0, \q \rho\in (0,2),
\]
for  large $x$. It is shown in \cite{YS} that for such potentials 
\[
\ln {\bf w}(\lambda) \sim -w_{0}\lambda^{-(2-\rho)/(2\rho)},\q w_{0}=w_{0}(b_{0},\rho)>0,
\]
as $\lambda\to + 0$. In this case condition \e{eq:Szeg1} is satisfied if and only if $\alpha< 3/2- 1/\rho$.  Observe that formula \e{eq:PolCou} for Pollaczek  polynomials corresponds to the exponential decay as $\lambda\to 0$ of the weight function ${\bf w}(\lambda)$ for Coulomb  repulsive potentials.

 The trace class condition for the operator ${\bf H}$ requires that $\rho>1$ which yields \e{eq:Szeg1} for any $\alpha<1/2$. Similarly, 
 the Hilbert-Schmidt condition for the operator ${\bf H}$ requires that $\rho>1/2$ which yields \e{eq:Szeg1} for any $\alpha< -1/2$. This is very close to the so-called quasi-Szeg\"o condition established in \cite{KS} for the Jacobi operator $H=H_{0}+V$ with a  Hilbert-Schmidt perturbation $V$.

  \section{Example. The Jacobi polynomials}
  
{\bf 8.1.}
In this section we suppose that the measure $d\rho (\lambda)$ is absolutely continuous and supported on the interval $[-1,1]$ so that
   \begin{equation}
d\rho (\lambda)= w(\lambda) d\lambda, \q  \lambda\in (-1,1),
\label{eq:Jac}\end{equation}
and
 \begin{equation}
 w(\lambda)= \kappa (1-\lambda)^{\alpha} (1+\lambda)^{\beta},  \q \alpha ,\beta>-1.
\label{eq:Jac1}\end{equation}
The weight function $ w(\lambda)= w_{\alpha ,\beta}(\lambda)$ (as well as all other objects discussed below) depends on $\alpha $ and $\beta$, but these parameters are often omitted in notation. 
 The constant 
 \begin{equation}
  \kappa= \kappa_{\alpha ,\beta}=\frac{\Gamma(\alpha+ \beta+2)}{2^{\alpha+ \beta+1}\Gamma(\alpha+1) \Gamma(\beta+1)}
\label{eq:Jac2}\end{equation}
 is chosen in such a way that the measure \e{eq:Jac} is normalized, i.e., 
$\rho ({\Bbb R})=\rho ((-1,1))=1$. 
The orthonormal polynomials $P_{n}(\lambda)= P_{n}^{(\alpha ,\beta)}(\lambda)$ constructed by the measure
\e{eq:Jac}, \e{eq:Jac1} are known as the Jacobi polynomials.

Usually   the Jacobi polynomials  ${\sf P}_{n} (\lambda) = {\sf P}_{n}^{(\alpha ,\beta)}(\lambda)$ are normalized (see, e.g., the books \cite{BE, Sz}) by the condition 
 \begin{equation}
  {\sf P}_{n} (1)=\frac{\Gamma(\alpha+1+n)}{ \Gamma(\alpha +1 ) \Gamma(n+1)}.
\label{eq:NNn}\end{equation}
According to formulas (10.8.4) and (10.8.5)
 in  the  book \cite{BE} we have 
  \begin{multline}
   \int_{-1}^{1}|{\sf P}_{n} (\lambda)|^{2}(1-\lambda)^{\alpha} (1+\lambda)^{\beta} d\lambda
   \\
   =\frac{2^{\alpha+\beta+1}\Gamma(\alpha+1+n)\Gamma(\beta+1+n)}{n! (2n+\alpha+\beta+1) \Gamma(\alpha +\beta+n+1 ) }=: {\sf h}_{n}
\label{eq:NN}\end{multline}
 and
   \[
{\sf P}_{n} (z)= {\sf k}_{n}
\big(z^{n}+ \frac{n(\alpha-\beta )}
{2 n +\alpha+\beta }z^{n-1}+\cdots\big)
\]
where
  \begin{equation}
{\sf k}_{n}= \frac{ \Gamma(\alpha+\beta+2n+1)}{2^n n! \Gamma(\alpha+\beta+n+1)}.
\label{eq:norm3}\end{equation}
For $\alpha+\beta\leq -1$ and $n=0$, we have to set
 \begin{equation}
{\sf h}_{0}= 2^{\alpha+\beta+1}\frac{ \Gamma(\alpha+1)\Gamma(\beta+1)}{  \Gamma(\alpha +\beta+2 ) }=\kappa^{-1} \q{\rm and}  
\q {\sf k}_0=1.
\label{eq:kh}\end{equation}
It follows from equalities \e{eq:Jac1}  and \e{eq:NN} that
   \[
P_{n} (z)=  (\kappa h_{n})^{-1/2}{\sf P}_{n} (z).
\]

Let $H$ be the Jacobi operator  corresponding to the measure \e{eq:Jac}, \e{eq:Jac1}. Its coefficients are given by formulas  \e{eq:norm1} which yields
   \begin{equation}
 a_{n}^{2} =  \frac{ {\sf k}_{n}^{2}}{ {\sf k}_{n+1}^{2}}   \frac{{\sf h}_{n+1}}{{\sf h}_{n }}
=\frac{4 (n+1)(n+\alpha+1)(n+\beta+1) }
 { (2n+\alpha+\beta+3) (2 n+\alpha+\beta  +2)^2}
 \,
  \frac{ n+\alpha+\beta+1}
 {2 n+\alpha+\beta+1 }
\label{eq:norm4}\end{equation}
(if $n=0$ and $\alpha+\beta=-1$ the last factor here should be replaced by $1$) and
   \begin{equation}
b_{n}  = (\alpha-\beta)\big(\frac{n}{2n+ \alpha+\beta} -\frac{n+1}{2n+ 2+\alpha+\beta}\big) 
\label{eq:norm5}\end{equation}
(if $n=0$ and $\alpha+\beta= 0$ this formula should be replaced by $b_{0}= (\beta-\alpha)/2$).
Obviously, for a symmetric weight function  \e{eq:Jac1} when $\alpha=\beta$, we have $b_{n} =0$.  
In the general case, we have the asymptotic formulas
 \begin{equation}
 a_{n}  =1/2 + 2^{-4} (1 -2\alpha^2-2\beta^2) n^{-2} +O\big( n^{-3}\big),\q
b_{n}  = 2^{-2}  (\beta^2-\alpha^2) n^{-2} +O\big( n^{-3}\big). 
\label{eq:norm6}\end{equation}
Assumption \e{eq:Tr}  is of course true in this case, and the sequence \e{eq:sig} satisfies the condition $\rho_{n} =O(n^{-1})$  as $n\to\infty$. Note that the  coefficients at $ n^{-2}$ in \e{eq:norm6} are not zeros unless $|\alpha|=|\beta|=1/2$.

Next, we calculate the product \e{eq:AA}.  It follows from \e{eq:norm4} that   
 \[
A^{(n)}:= \prod_{k=0}^{n-1} (2a_{k})=   \frac{ {\sf k}_{0}}{ \sqrt{{\sf h}_{0}}}   \,  \frac { \sqrt{{\sf h}_{n}}} { {\sf k}_{n}} \, 2^n.
\]
Let us find the limit of  the right-hand side as $n\to\infty$. According to formula (1.18.4) in the book~\cite{BE}, for all $p\in{\Bbb R}$ we  have
 \begin{equation}
 \Gamma (n+p)=n^p\, \Gamma( n)(1+O(n^{-1})).
\label{eq:AA2}\end{equation}
Therefore the sequence ${\sf h}_n$ defined by \e{eq:NN}  satisfies
 \begin{equation}
{\sf h}_n= 2^{\alpha+\beta} n^{-1}  + O(n^{-2}).
\label{eq:AA2x}\end{equation}
Let us use  additionally (see formula (1.3.11) in  \cite{BE}) that
\[
\Gamma(2n+\alpha+\beta+ 1)=\pi^{-1/2}2^{2n+\alpha+\beta}\Gamma(n+ (\alpha+\beta+ 1)/2)\Gamma(n+ (\alpha+\beta+ 2)/2).
\]
So, we see that the sequence ${\sf k}_n$ defined by \e{eq:norm3}  satisfies
$
{\sf k}_n=2^{\alpha+\beta}  2^n (\pi n)^{-1/2}   (1+ O(n^{-1}))$.
It follows that
 \[
\lim_{n\to\infty}2^n \frac { \sqrt{{\sf h}_n}} { {\sf k}_n} =\sqrt{\pi}2^{-(\alpha+\beta)/2}.
\]
Taking  also into account  formulas \e{eq:kh}, we finally find that
 \begin{equation}
A =\lim_{n\to\infty} A^{(n)}=\sqrt{\pi\kappa}\, 2^{-(\alpha+\beta)/2}
\label{eq:AB}\end{equation}
where $\kappa$ is defined by \e{eq:Jac2}.

\bigskip

{\bf 8.2.}
Let us   calculate the Szeg\"o function defined by \e{eq:SzY} for the measure  \e{eq:Jac},  \e{eq:Jac1}.  Now we have
 \[
w (\cos \theta)|sin \theta|= \kappa\, 2^{\alpha+\beta+1}|\sin(\theta/2)|^{2\alpha+1} (\cos(\theta/2))^{2\beta+1}.
\]
Substituting this expression into \e{eq:SzY} and using relations \e{eq:wwz},    we find that
  \begin{equation}
D (\z) =\sqrt{\kappa} 2^{-(\alpha+\beta+1)/ 2} (1-\z  )^{\alpha+1/2}   (1+\z  )^{\beta +1/2 }.
\label{eq:SAD}\end{equation}
In particular, for the operator $H_{0}$ when $\alpha=\beta=1/2$, we have
  $D_{0}(\z) = (2\pi)^{-1/2} (1-\z^{2})$. Setting in \e{eq:SAD} $\z=0$ and taking into account that, by \e{eq:SzD3},
   $\sqrt{2\pi} D(0) =A$, we recover expression \e{eq:AB}.

Now the perturbation determinant \e{eq:Tr1}  for the pair $H_{0}$, $H$   can be found from  formulas \e{eq:SzD3} where 
$B(\z)=1$ and
\e{eq:AB}, \e{eq:SAD}.
 Let us state the   results obtained.
    
  \begin{theorem}\label{CAJ}
 Let $H$ be the Jacobi operator with the coefficients \e{eq:norm4}, \e{eq:norm5} or, equivalently, corresponding to the spectral measure  \e{eq:Jac}, \e{eq:Jac1}. Then the Szeg\"o function $D (\z)$  and the perturbation determinant $\Delta (z) $ are  given by formulas \e{eq:SAD} and
  \begin{equation}
\Delta (z) = (1-\z (z))^{-\alpha+1/2}   (1+\z (z))^{-\beta +1/2 }.
\label{eq:AD}\end{equation}
 \end{theorem}

 Recall that the Jost function $ \Omega(z)$ is related to $\Delta(z)$ by formula \e{eq:HH}.  So it   follows from Theorem~\ref{GE1}   that
    \[
\lim_{n\to\infty} \z(z)^n P_{n}(z)= \kappa^{-1/2} \pi^{-1/2} 2^{(\alpha+\beta)/ 2} 
(1-\z (z))^{-\alpha-1/2}   (1+\z (z))^{-\beta -1/2 },\q z\in \Pi. 
\]
Thus we recover the well known asymptotic formula (see, for example, formula (8.21.9) in the book \cite{ Sz}) for the Jacobi polynomials.

  Let us now calculate $ \Omega (\lambda+i0)$ for $\lambda\in (-1,1)$; as usual, we set $\lambda=\cos\theta$ where $\theta\in (0,\pi)$. Then $\z=\z  (\lambda+i0) = e^{- i \theta} $   belongs to the lower half-circle. It is easy to see that
  \[
  1-\z= 1-\cos\theta + i\sin\theta=2 \sin(\theta/2) e^{i(\pi-\theta)/2}=\sqrt{2(1-\lambda)}  e^{i(\pi-\theta)/2}
  \]
  and
  \[
  1+\z= 1+\cos\theta - i\sin\theta=2 \cos(\theta/2) e^{-i\theta/2}=\sqrt{2(1+\lambda)} e^{-i\theta/2}.
  \]
Substituting these expressions into formula  \e{eq:AD} we obtain the following result. 

\begin{theorem}\label{CAJ1}
Under the assumptions of Theorem~\ref{CAJ}, the perturbation determinant  is given on the cut along $[-1,1]$
 by the  formula   
 \begin{equation}
\Delta(\lambda+i0)=   2^{(1-\alpha-\beta)/2}   (1-\lambda)^{(1-2\alpha)/4}(1+\lambda)^{(1-2\beta)/4} e^{i\pi \xi(\lambda)},\q \lambda\in (-1,1) ,
\label{eq:SZ+}\end{equation}
where the spectral shift function 
\begin{equation}
\xi(\lambda) =  (2\pi)^{-1} (\alpha+\beta-1) \arccos\lambda   - 4^{-1} (2 \alpha-1)  .
\label{eq:SZsf}\end{equation}
 \end{theorem}

 Using  \e{eq:HH}, \e{eq:AB} and substituting expression \e{eq:SZ+} for $|\Delta(\lambda+i0)|$ into formula \e{eq:SF1}, we recover relations 
\e{eq:Jac1}, \e{eq:Jac2}  for $w(\lambda)$.  
   In view of definition \e{eq:AP}, Theorem~\ref{CAJ1} yields     expressions for the limit amplitude  and   the limit phase.
   Thus,   formula \e{eq:Sz} means that
 \begin{multline*}
 P_{n} (\lambda)=  2^{1/2}( \pi\kappa)^{-1/2}(1-\lambda)^{-(1+2\alpha)/4}(1+\lambda)^{-(1+2\beta)/4} 
 \\
 \times \sin \big((n+\frac{\alpha+\beta+1}{2}) \arccos\lambda -\frac{2\alpha-1}{4} \pi\big) + O (n^{-1})
 \end{multline*}
where  $\kappa$ is given by  \e{eq:Jac2}.
This coincides of course with  asymptotics (8.21.10) of the Jacobi polynomials in the book~\cite{Sz}.
   
   Substituting expressions \e{eq:AD} for the perturbation determinant and \e{eq:SZsf} for the spectral shift function into formula  \e{eq:DD} we obtain a  curious identity
\[
\int_{-1}^1 \frac{\arccos \lambda}{\lambda-z}d\lambda=-2\pi \ln (1+\z (z)).
\]

The asymptotics of Jacobi polynomials at the edge points $z=\pm 1$ is different from general results  \e{eq:P+}  and  \e{eq:Vx1} of Section~4 because in view of formulas  \e{eq:norm5}  and  \e{eq:norm6} the condition \e{eq:Trfd}  is not satisfied unless $|\alpha|=|\beta|=1/2$. Indeed, putting together formulas \e{eq:NNn}, \e{eq:AA2} and \e{eq:AA2x}, we see that
 \begin{equation}
 P_{n} (1)= \kappa^{-1/2} \Gamma(\alpha+1)^{-1}2^{-(\alpha+\beta)/2} n^{ \alpha+1/2} \big(1+O(n^{-1})\big)
\label{eq:th}\end{equation}
where $\kappa$ is defined by \e{eq:Jac2}. This is consistent with asymptotics \e{eq:P+} if $\alpha=1/2$ and with asymptotics \e{eq:Vx1} if $\alpha=-1/2$. The results for 
$ P_{n} (- 1)$ are quite
 similar.

\bigskip

{\bf 8.3.}
Let us now discuss stationary representations for the wave operators $W_{\pm} (H,H_{0})$ and the corresponding 
scattering matrix  $S (\lambda)$. Recall that the  operator $F$ was defined  by equations  \e{eq:UF} and \e{eq:psX}.     Using formulas \e{eq:CH1r} and  \e{eq:SS}, we can state the following result.

\begin{theorem}\label{WOS}
Under the assumptions of Theorem~\ref{CAJ},
  the wave operator $W_{\pm} (H,H_{0})$ is given by the equality   \e{eq:CH1pr} where $\Sigma_{\pm}$ is the operator of multiplication by the function
\[
\sigma_{\pm} (\lambda)=e^{\pm i (   (\alpha+\beta-1) \arccos\lambda -\pi ( \alpha-1/2) )/2}.
\]
The scattering matrix  satisfies the equation
\[
S (\lambda)=e^{i\pi ( \alpha-1/2)  } e^{- i    (\alpha+\beta-1) \arccos\lambda}  .
\]
 \end{theorem}

In view of asymptotics \e{eq:norm6} the condition \e{eq:Trfd} is not satisfied (unless $|\alpha|=|\beta|=1/2$), and hence  the results of Subsections~4.1 and 4.2, as well as Theorem~\ref{Levinson}, are not applicable. Nevertheless representation  \e{eq:AD} allows us to find the singularity of the perturbation determinant $\Delta(z) $ as $z\to\pm 1$.  If $z\to  1$, then
\[
\Delta(z)= 2^{-\beta+1/2} (1-\z(z))^{-\alpha +1/2}\big( 1+ O (|\z(z)-1|)\big).
\] 
If $z\to - 1$, we have
\[
\Delta(z)= 2^{-\alpha+1/2} (1+\z(z))^{-\beta +1/2}\big( 1+ O (|\z(z)+1|)\big).
\] 
 It follows  from \e{eq:SZsf} that the limits of   $\xi(\lambda)$ as $\lambda\to -1+0$ and as $\lambda\to 1-0$ exist and
\[
\xi (-1+0)=(1-2\beta)/4,   \q \xi (1-0)= (1-2 \alpha)/4.
\]
Recall that $\xi(\lambda)=0$ for  $|\lambda|>1$. Thus, the spectral shift function is continuous at the point $-1$ (at the point $1$) if and only if $\beta=1/2$ ($\alpha=1/2$).  
We are not aware of the     results of this type for the differential operator $D a(x) D + b(x)$ in the space $L^2 ({\Bbb R}_{+})$.  

Let us, finally, discuss the exceptional case $|\alpha|=|\beta|=1/2$:

$1^{0}$ If $\alpha=\beta=1/2$, then $ a_{n}  =1/2$, $b_{n}=0$   for all $n\in{\Bbb Z}_{+}$, and hence $H=H_{0}$.

 $2^{0}$ If $\alpha=\beta=-1/2$, then $a_{0}=1/\sqrt{2}$, $ a_{n}  =1/2$  for all $n\geq 1$,  $b_{n}=0$   for all $n\in{\Bbb Z}_{+}$ and the perturbation $V$ has rank $2$. In this case, $\Delta(z)=1-\z(z)^{2}$ so that the corresponding Jacobi operator has resonances at both points $z=1$ and $z=-1$.
 
 $3^{0}$ If $\alpha=-\beta=\pm 1/2$, then $ a_{n}  =1/2$   for all $n\in{\Bbb Z}_{+}$, $b_{0}=\mp 1$, 
$ b_{n}  =0$   for all $n\geq 1$ and the perturbation $V$ has rank $1$. In this case, $\Delta(z)=1\pm \z(z)$ so that the corresponding Jacobi operator has a resonance at the  point $z=\mp 1$.

Observe  that
in   cases $1^{0}$ and $2^{0}$, $P_{n} (\lambda)$ are   the  normalized  Chebyshev polynomials of the second and first kind, respectively.

\bigskip

 {\bf 8.4.}
 Note that very detailed results on the asymptotics of the orthogonal polynomials were obtained in \cite{McL} for the weight function 
   \[
 \tilde{w}(\lambda)=   (1-\lambda)^{\alpha} (1+\lambda)^{\beta}h(\lambda);
\]
here $h(\lambda)$ is an arbitrary real analytic function such that  $h(\lambda)>0$ for $\lambda\in[-1,1]$. On the contrary, it was shown in \cite{M-F} that, for the weight function $ \tilde{w}(\lambda)\chi_{c}(\lambda)$ where $\chi_{c}(\lambda)=1$ for $\lambda\in [-1,0)$ and $\chi_{c}(\lambda)=c\neq 1$, $c>0$,  for $\lambda\in (0,1]$,  the classical asymptotics \e{eq:zF} of the orthogonal polynomials is significantly changed. For such weight functions,  $ a_{n}-1/2 =O (n^{-1})$ and $b_{n} =O (n^{-1})$ so that condition \e{eq:Tr} is not satisfied.

 Both papers \cite{McL}  and \cite{M-F}  rely on the approach of \cite{F-I-K, D-Z}.

 \appendix
 
 \section{Discrete Volterra equations}

{\bf A.1.}
Here  we prove Theorems~\ref{Jost} and \ref{JostP}. We follow rather closely the scheme exposed for the Schr\"odinger equation in a detailed way in Section~4.1  of \cite{YA}. Note however that the operator $H$ is obtained by a   second order perturbation of the operator $H_{0}$. This circumstance should be taken into account.
  
  Let us reduce equation \e{eq:Jy} to a discrete Volterra integral equation. 
    
  \begin{lemma}\label{Jos1}
Let assumption   \e{eq:Tr} be satisfied, and let $z\in\clos{\Pi}$, $z\neq \pm 1$. Suppose that $f(z)=\{f_{n} (z)\}\in \ell^\infty ({\Bbb Z}_{+})$ and     define $V$ by formula \e{eq:JV}. Then a solution of the equation \e{eq:V}
satisfies also equation \e{eq:Jy}.
\end{lemma}

 \begin{pf}
 Under assumption   \e{eq:Tr}, the sequence $ \phi(z) : =V f(z) \in \ell^1 ({\Bbb Z}_{+})$.
Put
  \[
   \Sigma^{(+)}_{n}(z)= \sum_{m=n+1}^\infty   \z(z)^{m } \phi_{m}(z),\q \Sigma^{(-)}_{n}(z)= \sum_{m=n+1}^\infty   \z(z)^{-m }  \phi_{m}(z).
\]
Then equation \e{eq:V} can be written as
\[
 f_{n} =\z ^n -\frac{\z^{n }}{\sqrt{ z^2-1}}\Sigma^{(-)}_{n} + \frac{\z^{-n }}{\sqrt{ z^2-1}}\Sigma^{(+)}_{n},
 \]
whence
 \begin{align*}
 f_{n-1}& =\z ^{n-1} -\frac{\z^{n-1} }{\sqrt{ z^2-1}}\Sigma^{(-)}_{n} + \frac{\z^{-n +1}}{\sqrt{ z^2-1}}\Sigma^{(+)}_{n} + \frac{\z-\z^{-1 }}{\sqrt{ z^2-1}}  \phi_{n},
 \\
 f_{n+1}& =\z ^{n+1} -\frac{\z^{n +1}}{\sqrt{ z^2-1}}\Sigma^{(-)}_{n} + \frac{\z^{-n-1 }}{\sqrt{ z^2-1}}\Sigma^{(+)}_{n}.
 \end{align*}
 It follows that
  \begin{multline}
2 \sqrt{ z^2-1} ((H_{0}-z) f)_{n}=\sqrt{ z^2-1} \big(f_{n-1}+f_{n+1}  - 2 z f_{n}\big)= (\z-\z^{-1})  \phi_{n} 
\\
+  (\z^{-n+1} +\z^{-n-1})\Sigma^{(+)}_{n} -(\z^{n+1} + \z^{n-1})\Sigma^{(-)}_{n} -2 z(\z^{-n}\Sigma^{(+)}_{n}-\z^{n}\Sigma^{(-)}_{n} ) .
\label{eq:AX1}\end{multline}
The sums of the terms containing $\Sigma^{(+)}_{n}$ and $\Sigma^{(-)}_{n}$ equal zero because $\z+\z^{-1}=2z$.  Since
$\z-\z^{-1}=-2\sqrt{ z^2-1}$, equality \e{eq:AX1} yields the desired equation   $((H_{0}-z) f)_{n}=-(V f)_{n}$.
 \end{pf}
 
Next, we study equation \e{eq:V}. Let us distinguish  the term corresponding to $m=n+1$ in   its right-hand side. Then equation    \e{eq:V} reads as
  \begin{multline*}
2 a_{n}  f_{n}(z)=\z(z)^n -2b_{n+1}  f_{n+1} (z) +   (1-2 a_{n+1}) f_{n+2} (z)
\\
 -\frac{1}{\sqrt{ z^2-1}}\sum_{m=n+2}^\infty (\z(z)^{n-m}- \z(z)^{m-n}) (Vf(z))_{m}.
  \end{multline*} 
It is convenient to rewrite this equation  in terms of the  sequence
\begin{equation}
  g_{n}(z)=2 a_{n}\z(z)^{-n}  f_{n}(z)
\label{eq:V1}\end{equation}
as
\begin{equation}
  g_{n}(z)=1 + \sum_{m=n +1}^\infty G_{n,m} (z) g_{m}(z)
\label{eq:V3}\end{equation}
where
\begin{multline}
 2a_{m}\sqrt{z^2-1}    G_{n,m} (z)= 
  (\z(z)^{2 m-2n}-1) b_{m} 
 \\
  +  ( \z(z)^{2 m-2n-1}-\z(z)) (a_{m-1}  -1/2) +(  \z(z)^{2 m-2n+1}-\z(z)^{-1})     (a_{m}  -1/2) .
 \label{eq:V2}\end{multline}
It is important that
  \begin{equation}
| G_{n,m} (z)|\leq C \alpha_{m}  
\label{eq:V4}\end{equation}
where 
  \[
 \alpha_{m}=|a_{m-1}-1/2| + |a_m-1/2| +  |b_m|
\]
and the constant $C=C(z)$   does not depend on $n$, $m$ and on  $z$ in compact subsets of $\clos\Pi$ away from the points $\pm 1$.

Let us solve equation \e{eq:V3} by iterations. 
 Put   $g^{(0)}_n=1$ for all $n\in {\Bbb Z}_{+}$ and 
  \begin{equation}
 g^{(k+1)}_{n}(z)= \sum_{m=n +1}^\infty G_{n,m} (z) g^{(k )}_m(z),\q k\geq 0.
\label{eq:V5}\end{equation}

  \begin{lemma}\label{Jos2}
The estimates
 \begin{equation}
| g^{(k )}_{n}(z)|\leq \frac{C^k}{k!}\big(\sum_{m=n+1 }^\infty \alpha_{m}\big)^k,\q \forall n\in{\Bbb Z}_{+},
\label{eq:V6}\end{equation}
are true for all $k\in{\Bbb Z}_{+}$. 
\end{lemma}

  \begin{pf}
  Suppose that \e{eq:V6} is satisfied for some $k\in{\Bbb Z}_{+}$. We have to check that
 the same estimate (with $k$ replaced by $k+1$ in the right-hand side) holds for $ g^{(k+1)}_{n}(z)$.  
  It  follows from  \e{eq:V4} and  \e{eq:V6} that
   \begin{equation}
| g^{(k +1)}_{n}(z)|\leq \frac{C^{k+1}}{k!}  \sum_{m=n +1}^\infty  \alpha_{m}\big( \sum_{p=m+1}^\infty \alpha_{p}\big)^k.
\label{eq:V7}\end{equation}
Observe that
 \[
(k+1)   \alpha_{m}\big( \sum_{p=m+1}^\infty \alpha_{p}\big)^k + \big( \sum_{p=m+1}^\infty \alpha_{p}\big)^{k+1}
\leq
\big( \sum_{p=m }^\infty \alpha_{p}\big)^{k+1},
\]
and hence, for all $N\in{\Bbb Z}_{+}$,
   \begin{multline*}
 (k+1)  \sum_{m=n +1}^N  \alpha_{m}\big( \sum_{p=m+1}^\infty \alpha_{p}\big)^k
 \\
 \leq 
 \sum_{m=n +1}^N \Big( \big( \sum_{p=m }^\infty \alpha_{p}\big)^{k+1}-  \big( \sum_{p=m+1}^\infty \alpha_{p}\big)^{k+1}\Big)\leq  \big( \sum_{p=n+1}^\infty \alpha_{p}\big)^{k+1}.
 \end{multline*}
Substituting this estimate into   \e{eq:V7}, we obtain \e{eq:V6} for $g^{(k +1)}_{n}(z)$.
    \end{pf}
    
    It is now easy to conclude the proof of Theorem~\ref{Jost}.
    Set
     \[
   g_{n} (z)=\sum_{k=0}^\infty g^{(k)}_{n}(z).
\]
Estimate \e{eq:V6} shows that this series converges and, by the recurrent definition \e{eq:V5}, the functions $ g_{n} (z)$ satisfy
equation \e{eq:V3} equivalent to \e{eq:V}. Since every function $g^{(k)} _{n}(z)$ is analytic in $z\in\Pi$ and is continuous up to the cut $[-1,1]$ (away from the points $\pm 1$), estimate \e{eq:V6} guarantees that $ g_{n}(z)$ and  hence $f_{n}(z)$ possess the same properties.

It also follows from \e{eq:V3} and \e{eq:V4} that
\[
 |   g_{n}(z)-1 | \leq C\sum_{m=n +1}^\infty \alpha_{m} 
\]
which in view of  \e{eq:V1} implies  \e{eq:Jost}.
   $ \q  \Box$

As far as Theorem~\ref{JostP}  is concerned, we note that all functions $g^{(k)}_{n}(z)$ are continuous with respect to the cut-off parameter $N\to\infty$, and therefore estimates obtained in the proof of Theorem~\ref{Jost} imply relation \e{eq:cut}.

  \bigskip
  
  {\bf A.2.}  
  Next, we discuss Theorem~\ref{OPP}.
     The first assertion  is quite similar to Lemma~\ref{Jos1}, and so its proof will be omitted.

 \begin{lemma}\label{OPR2}
Let   $P(z)=\{P_{n} (z)\}_{n=-1}^\infty $ where $P_{-1} (z)=0$, $P_0 (z)=1$ satisfy equations \e{eq:OPR1}. Then    $P_{n}(z)$
satisfy also equations \e{eq:Py}; in particular,  $P_{n} (z)$ is a polynomial of degree $n$.
\end{lemma}

We proceed from   equation \e{eq:OPR1} and,
following the scheme of the previous subsection,   distinguish the term corresponding to $m=n-1$ in the sum \e{eq:OPR1}. Then equation \e{eq:OPR1} reads as
   \begin{multline*}
2 a_{n-1} P_{n}(z)=P_{n}^{(0)}(z) -2 \alpha_{n-2}P_{n-2}(z)  -2 b_{n-1}P_{n-1}(z)
\\
 +\frac{1}{\sqrt{ z^2-1}}\sum_{m=0}^{ n-2} (\z(z)^{n-m}- \z(z)^{m-n}) (VP(z))_{m},\q n \geq 2.
\end{multline*}
It is  convenient to rewrite this equation in terms of the sequence
   \[
 Q_{n}(z)= 2 a_{n-1} \z(z)^{n}  P_{n}(z)
\]
   as
\begin{equation}
 Q_{n}(z)=Q_{n}^{(0)}(z) - \sum_{m=0}^{n-1}  {\sf G}_{n,m} (z)Q_{m}(z).
\label{eq:M3}\end{equation}
Here $ Q_{n}^{(0)}(z)= \z(z)^{n}  P_{n}^{(0)}(z)$ 
are bounded uniformly in $n$ in view of \e{eq:OPR}. The matrix elements ${\sf G}_{n,m} (z)$ are quite similar to $G_{n,m} (z)$ defined by   \e{eq:V2}. It is only important that ${\sf G}_{n,m} (z)$   satisfy estimate \e{eq:V4}.

We again solve equation \e{eq:M3} by iterations. 
 Put      $Q^{(k )}_0 =0$ for $k\geq 1$ and
   \begin{equation}
 Q^{(k+1)}_{n}(z)= \sum_{m=0}^{n-1}  {\sf G}_{n,m} (z) Q^{(k )}_m(z),\q k\geq 0, \q n\geq 1.
\label{eq:M5}\end{equation}

  \begin{lemma}\label{M2}
  The estimates
 \begin{equation}
| Q^{(k )}_{n}(z)|\leq \frac{C^k}{k!}\big(\sum_{m=0}^{n-1}  \alpha_{m}\big)^k,\q \forall n \geq 1,
\label{eq:M6}\end{equation}
are true for all $k\in{\Bbb Z}_{+}$. 
\end{lemma}

  \begin{pf}
    Suppose that \e{eq:M6} is satisfied for some $k\in{\Bbb Z}_{+}$. We have to check that
 the same estimate (with $k$ replaced by $k+1$ in the right-hand side) holds for $ Q^{(k+1)}_{n}(z)$.  
  It  follows from  \e{eq:M5},  \e{eq:M6} that
   \begin{equation}
| Q^{(k +1)}_{n}(z)|\leq \frac{C^{k+1}}{k!}  \sum_{m=0}^{n -1}  \alpha_{m}\big( \sum_{p=0}^{m-1} \alpha_{p}\big)^k.
\label{eq:M7}\end{equation}
Observe that
 \[
(k+1)   \alpha_{m} \big( \sum_{p=0}^{m-1} \alpha_{p}\big)^k + \big( \sum_{p=0}^{m-1} \alpha_{p}\big)^{k+1}
\leq
\big( \sum_{p=0 }^m \alpha_{p}\big)^{k+1},
\]
and hence
\[
 (k+1)  \sum_{m=0}^{n-1}  \alpha_{m}\big( \sum_{p=0}^{m-1} \alpha_{p}\big)^k\leq 
 \sum_{m=0}^{n-1} \Big( \big( \sum_{p=0 }^m \alpha_{p}\big)^{k+1}-  \big( \sum_{p=0}^{m-1} \alpha_{p}\big)^{k+1}\Big)\leq  \big( \sum_{p=0}^{n-1} \alpha_{p}\big)^{k+1}.
\]
Substituting this estimate into   \e{eq:M7}, we obtain \e{eq:M6} for $k+1$.
    \end{pf}

      Let us now set
      \[
   Q_{n} (z)=\sum_{k=0}^\infty Q^{(k)}_{n}(z).
\]
Estimate \e{eq:M6} shows that this series converges and
\[
   Q_{n} (z)\leq\sum_{k=0}^\infty \frac{C^k}{k!}\big(\sum_{m=0}^{n-1}  \alpha_{m}\big)^k=
    \exp\big( C \sum_{m=0}^{n-1} \alpha_{m}\big).
\]
By the recurrent definition \e{eq:M5}, the functions $ Q_{n} (z)$ satisfy
equation \e{eq:M3} equivalent to \e{eq:OPR1}.
This proves estimate \e{eq:aa2}. $\q\Box$

   \bigskip
   
   {\bf A.3.} 
The proof of Theorem~\ref{JOSF} can be obtained similarly to that of Theorem~\ref{Jost} if one observes that
 \begin{equation}
\Big|\frac{ \z^{2k} -1}{\sqrt{z^2-1}} \Big|= 2 \Big|\z\frac{   \z^{2k} -1}{ \z^2-1} \Big|
= 2  |\z| |\z^{2(k-1)}+\z^{2(k-2)}+\cdots+\z^2+ 1|\leq 2 k.
\label{eq:ZR1}\end{equation}
Thus instead of \e{eq:V4}, we now have a bound 
 \begin{equation}
| G_{n,m} (z)|\leq C(m+1) \alpha_{m}  
\label{eq:ZR2}\end{equation}
with a constant $C$ not depending on $\z$ in the unit disc. Under assumption \e{eq:Trfd} this bound allows us to repeat the arguments of Theorem~\ref{Jost}. 
 
  The proof of Theorem~\ref{OPP1} follows the scheme of  proof of Theorem~\ref{OPP}, but now we again have to use estimate \e{eq:ZR1}. This estimate implies that the kernels ${\sf G}_{n,m} (z)$ in \e{eq:M3} obey bound \e{eq:ZR2}.


\begin{thebibliography}{99}
 
 \bibitem{AKH} N.~Akhiezer,
\emph{The classical moment problem and some related questions in analysis},  Oliver and  Boyd, Edinburgh and London, 1965.

\bibitem {Ber} Yu. M. Berezanskii,  {\it Expansion in eigenfunctions of selfadjoint operators}, Amer. Math. Soc., Providence, R.I., 1968.


 \bibitem {Bern} S. Bernstein,
  {\em Sur les polyn\^omes orthogonaux relatifs  \`a un segment fini}, Journal de Math\'ematiques, {\bf 9}  (1930), 127-177; 
  {\bf 10}  (1931), 219-286.



 \bibitem {BF}V. S. Buslaev and L. D. Faddeev,
  {\em Formulas for traces for a singular Sturm-Liouville differential operator}, Soviet Math. Dokl., {\bf 1}  (1960), 451-454.

 \bibitem {Case} K.~M.~Case,   {\em Orthogonal polynomials.} II, 
 Journal of Math. Phys., {\bf 16} (1975),  1435-1440.

 

\bibitem {D-K} D.~Damanik,  R.~Killip, {\em Half-line Schr\"odinger operators with no bound states}, Acta Math., {\bf 193}  (2004), 31-72. 


   \bibitem{D-Z}     P.~Deift, X.~Zhou, {\em A steepest descent method for oscillatory Riemann-Hilbert problem},  Ann.
Math., {\bf 137} (1993), 295-368.


    \bibitem{Duren} P.~L.~Duren,  {\em Theory of $H^p$ spaces}, Academic Press, New York and London, 1970.
 
 \bibitem{BE} A.~Erd\'elyi,  W.~Magnus, F.~Oberhettinger, F.~G.~Tricomi,
\emph{Higher transcendental functions}, Vol. 1, 2, 
 McGraw-Hill, New York-Toronto-London, 1953.
 
  \bibitem {Fadd1} L. D. Faddeev, {\em An expression for the trace of the difference between two singular differential operators of Sturm-Liouville type}, Dokl. AN SSSR    {\bf 115},  N 5    (1957),  878-881 (Russian).
 
 

\bibitem {F1} L. D. Faddeev,  {\em Properties of the $S$-matrix of the one-dimensional Schr\"odinger equation}, Amer. Math. Soc. Transl. (Ser. 2)  {\bf 65}  (1967), 139-166.

\bibitem {Finv} L. D. Faddeev,  {\em Inverse problem of quantum scattering theory.  II}, J. Soviet. Math., {\bf 5} (1976), 334-396.

\bibitem {F-T} L.~Faddeev, L.~Takhtajan, {\em Hamiltonian methods in the theory of solitons}, Springer-Verlag, Berlin, Heidelberg, 2007.
 

\bibitem {F-I-K} A. Fokas, A. Its, A. Kitaev,   {\em The isomonodromy approach to matrix models in $2D$ quantum gravity}, 
Comm.
Math. Phys., {\bf 147} (1992), 395-430. 

\bibitem {Ge-Le} I. M. Gel'fand, B. M. Levitan, {\em  On the determination of a differential equation from its spectral function}, Izv. Akad. Nauk SSSR, Ser. Mat. {\bf 15}  (1951), 309-361; Amer. Math. Soc. Transl. (Ser. 2)  {\bf 1}  (1955), 253-304.


\bibitem {GK} I. C. Gokhberg and M. G. Kre\u{\i}n, {\em Introduction to the theory of linear
nonselfadjoint operators in Hilbert space}, Amer. Math. Soc., Providence, Rhode Island, 1970.

\bibitem {Gus} G. Sh. Guseinov,  {\em The scattering problem for an infinite Jacobi matrix}, Izv. Akad. Nauk Armyan. SSR Ser. Mat. {\bf 12} (1977), 365-379 (Russian).




 
 
       \bibitem{Hu-S}     D.~Hundertmark, B.~Simon, {\em Lieb-Thirring inequalities for Jacobi matrices},  J.
Approx. Theory, {\bf 118} (2002), 106-130.

    \bibitem {J-P} P. Jost, A. Pais, {\em  On the scattering of a particle by a static potential}, Phys. Rev., {\bf 82}  (1951), 840-851.
    
 

  \bibitem{KS}     R.~Killip, B.~Simon, {\em Sum rules for Jacobi matrices and their applications to spectral theory},  Ann.
Math., {\bf 158} (2003), 253-321.

  \bibitem{McL}     A. B. J.~Kuijlaars, K. T.-R.~McLaughlin, W. Van Assche, M. Vanlessen, {\em The Riemann-Hilbert approach to strong asymptotics for orthogonal polynomials on $[-1,1]$},  Adv.
Math. {\bf 188} (2004), 337-398.

\bibitem {LL} L. D. Landau and E. M. Lifshitz, {\it Quantum mechanics}, Pergamon Press, 1965.

  \bibitem{M-F} F.~Moreno, A.~Martinez-Finkelshtein, V.~L.~Sousa, {\em Asymptotics of orthogonal polynomials for a weight with a jump on $[-1,1]$}, Constr. Approx. {\bf 33}, No. 2,  (2011),  219-263. 

 \bibitem{Ryck} E.~Ryckman, {\em  A strong Szeg\"o theorem for Jacobi matrices},  Comm. Math. Phys. {\bf  271}  (2007),
791-820;   Erratum,  Comm.  Math. Phys. {\bf  275}  (2007),
581-585.
  
  \bibitem {Sim}B. Simon, {\em Trace ideal methods}, London Math. Soc. Lecture Notes, Cambridge Univ.
Press,  London and New York, 1979.


 
    \bibitem{Sz} G.~Szeg\"o,  {\em Orthogonal polynomials}, Amer. Math. Soc., Providence, R. I., 1978.
 
 \bibitem {Teschl}  G.~Teschl, {\em Jacobi operators and completely integrable nonlinear lattices}, Amer. Math. Soc.,   Providence,  R. I., 2000.
    


 

 


 \bibitem {YS}D. R. Yafaev, {\em The low energy scattering for slowly decreasing 
  potentials},  Comm.
Math. Phys. {\bf 85} (1982), 177-196. 


 \bibitem{Ya} D. R. Yafaev, {\em Mathematical scattering theory: General theory}, Amer. Math. Soc.,   Providence,
  R. I., 1992.
  
   \bibitem{YA} D. R. Yafaev, {\em Mathematical scattering theory: Analytic  theory}, Amer. Math. Soc.,   Providence,
  R. I., 2010.
  

 \bibitem{Yjmp} D. R. Yafaev, {\em A point interaction for the discrete Schr\"odinger operator and generalized Chebyshev polynomials}, Journal of Math. Phys., {\bf 58} (2017),  063511.
    
 

     
    \end{thebibliography}
\end{document}